\documentclass[a4paper,11pt]{article}
\usepackage{latexsym,amssymb}
\usepackage{amsmath}
\usepackage{amsthm}
\usepackage{cases}
\usepackage[mathscr]{eucal}
\usepackage{color}
\usepackage[abbrev]{amsrefs}

\newtheorem{Theorem}{\bf Theorem}[section]
\newtheorem{Lemma}{\bf Lemma}[section]
\newtheorem{Proposition}{\bf Proposition}[section]
\newtheorem{Corollary}{\bf Corollary}[section]
\newtheorem{Remark}{\bf Remark}[section]
\newtheorem{Example}{\bf Example}[section]
\newtheorem{Definition}{\bf Definition}[section]

\newenvironment{theorem}{\begin{Theorem}$\!\!\!$}{\end{Theorem}}
\newenvironment{lemma}{\begin{Lemma}$\!\!\!$}{\end{Lemma}}

\newenvironment{corollary}{\begin{Corollary}$\!\!\!$}{\end{Corollary}}
\newenvironment{remark}{\begin{Remark}$\!\!\!$}{\end{Remark}}

\newenvironment{definition}{\begin{Definition}$\!\!\!$}{\end{Definition}}

\renewcommand{\qed}{{\hfill$\Box$ }}

\newcommand{\N}{\mathbb{N}}
\newcommand{\R}{\mathbb{R}}

\numberwithin{equation}{section}

\def\Xint#1{\mathchoice
{\XXint\displaystyle\textstyle{#1}}%
{\XXint\textstyle\scriptstyle{#1}}%
{\XXint\scriptstyle\scriptscriptstyle{#1}}%
{\XXint\scriptscriptstyle\scriptscriptstyle{#1}}%
\!\int}
\def\XXint#1#2#3{{\setbox0=\hbox{$#1{#2#3}{\int}$}
\vcenter{\hbox{$#2#3$}}\kern-.5\wd0}}

\def\dashint{\Xint-}

\title{Solvability of Superlinear Fractional  Parabolic Equations}
\author{\qquad\\
Yohei Fujishima, Kotaro Hisa, Kazuhiro Ishige, and Robert Laister}
\date{}
\begin{document}
\maketitle
\begin{abstract}
We study necessary conditions and sufficient conditions
for the existence of local-in-time solutions of the Cauchy problem
for superlinear fractional parabolic equations.
Our conditions are sharp and clarify the relationship between
the solvability of the Cauchy problem and the strength of the singularities of the initial measure.
\end{abstract}
\vspace{20pt}
\smallskip
\noindent Y. F.:  Department of Mathematical and Systems Engineering, Faculty of Engineering,\\
\qquad\quad Shizuoka University, 3-5-1, Johoku, Hamamatsu, 432-8561, Japan.\\
\noindent
E-mail: {\tt fujishima@shizuoka.ac.jp}\\

\smallskip
\noindent K. H.:  Graduate School of Mathematical Sciences, The University of Tokyo,\\
\qquad\quad 3-8-1, Komaba, Meguro-ku, Tokyo 153-8914, Japan.\\
\noindent
E-mail: {\tt hisak@ms.u-tokyo.ac.jp}\\

\smallskip
\noindent K. I.:  \,Graduate School of Mathematical Sciences, The University of Tokyo,\\
\qquad\quad 3-8-1, Komaba, Meguro-ku, Tokyo 153-8914, Japan.\\
\noindent
E-mail: {\tt ishige@ms.u-tokyo.ac.jp}\\

\smallskip
\noindent 
R. L.: Department of Computer Science and Creative Technologies,\\
 \qquad\quad University of the West of England, Bristol BS16 1QY, UK.\\
\noindent
E-mail: {\tt Robert.Laister@uwe.ac.uk}\\
\vspace{15pt}
\vspace{3pt}
\newpage
\section{Introduction}
We consider the Cauchy problem for a  superlinear fractional  parabolic equation 
\begin{equation}
\tag{P}
\left\{
\begin{array}{ll}
\partial_t u+(-\Delta)^{\frac{\theta}{2}} u=F(u),\quad & x\in{{  \R}}^N,\,\,\,t>0,\vspace{3pt}\\
u(0)=\mu &  \mbox{in} \quad {{  \R}}^N,\vspace{3pt}
\end{array}
\right.
\end{equation}
where $\mu$ is a nonnegative Radon measure in ${{  \R}}^N$. 
Throughout the paper we assume that $N\ge 1$, $0<\theta\le 2$, and 
$F:[0,\infty)\to [0,\infty)$  is 
(at least) continuous. 

In general, the existence of local-in-time nonnegative solutions of  problem~(P)  depends crucially on the   delicate interplay
 between the  strength of the singularities of the initial measure $\mu$ and the behavior of $F(\tau)$ as $\tau\to\infty$.
In this paper, for a large class of 
 nonlinearities $F$, we obtain new  necessary conditions and new sufficient conditions
for the local solvability of  problem~(P). 
The prototypical example we have in mind is
$$
\mbox{$F(\tau)=\tau^p[\log(L+\tau)]^q$, 
     where $p>1$, $q\in{{  \R}}$, and $L\ge 1$}. 
$$
As a consequence of our more general results, we are then able to derive sharp results
for classes of nonlinearities which include these prototypes as special cases, and quantify this interplay more precisely via  `optimal singularities'. 

Throughout this paper we use the following notations. 
For $T>0$ we  set  $Q_T:={{\R}}^N\times (0,T)$ and  let $B(x,\sigma)$ denote the Euclidean ball in ${{  \R}}^N$ centre $x$, radius $\sigma$. We use $\dashint_{B}f\, dx$ for the average value of $f$ over $B$ with respect to  the Lebesgue measure  $dx$. 
The set of nonnegative Lebesgue measurable functions in ${\mathbb R}^N$ is denoted by ${\mathcal L}_0$,   while 
${\mathcal M}$ denotes the set of nonnegative Radon measures in ${\mathbb R}^N$. For $\mu\in {\mathcal L}_0$ we abuse
 terminology somewhat by speaking of   `measure~$\mu$'  defined via $d\mu =\mu (x)dx$.
\subsection{Background}
\label{subsection:1.1}
The solvability of the Cauchy problem for superlinear parabolic equations
has been studied in many papers since the pioneering work by Fujita~\cite{F}. 
The literature is now very extensive and we  refer to the comprehensive monograph~\cite{QS}. 
We also mention the following works, some of which are directly related to  this paper, others with a different emphasis (higher order equations, systems, nonlinear boundary conditions): superlinear  parabolic equations \cites{BP, BC, BF, F, KY, LRSV, LS, M, RS, T, Wang, W1, W2};  
  linear heat equation with  nonlinear boundary conditions \cites{DFL, GL, HI02, IS01, IS02} ;
superlinear parabolic equations with a potential  \cites{BK, BTW, C, HS, HT, Wang};
   superlinear   parabolic systems \cites{EH, FI01, FI02, IKS, QS01};
  superlinear fractional parabolic equations \cites{GK, HI01, HIT, LS02, S};
  superlinear higher order parabolic equations \cites{CM, GP, GG, IKO, IKO2}.

In \cite{HI01} the second and  third authors of this paper
 considered problem~(P) in the special case of the power law nonlinearity $F(u)=u^p$ with $p>1$:
\begin{equation}
\label{eq:1.1}
\left\{
\begin{array}{ll}
\partial_t u+(-\Delta)^{\frac{\theta}{2}} u=u^p,\quad & x\in{{  \R}}^N,\,\,\,t>0,\vspace{3pt}\\
u(0)=\mu &  \mbox{in} \quad {{  \R}}^N.\vspace{3pt}
\end{array}
\right.
\end{equation}
There, as here, 
the exponent $p_\theta:=1+\theta/N$ plays a critical role.  
They proved the following necessary conditions for the local existence (cases (i) and (ii)). 
 \begin{itemize} 
  \item[(i)] 
  Let  $\mu\in{\mathcal M}$. 
  If problem~\eqref{eq:1.1} possesses a nonnegative solution in $Q_T$ for some $T>0$, 
  then there exists 
  $C_1=C_1(N,\theta,p)>0$ such that
  \begin{equation}
  \label{eq:1.2}
  \sup_{x\in{{  \R}}^N}\mu (B(x,\sigma))\le C_1\sigma^{N-\frac{\theta}{p-1}},
  \quad
  0<\sigma\le T^{\frac{1}{\theta}}.
  \end{equation}
 In the case where $1<p<p_\theta$
  the function $(0,\infty)\ni\sigma\mapsto \sigma^{N-\theta/(p-1)}$ is decreasing so that
 relation~\eqref{eq:1.2} is equivalent to
 $$
 \sup_{x\in{{  \R}}^N}\mu(B(x,T^{\frac{1}{\theta}}))\le C_1 T^{\frac{N}{\theta}-\frac{1}{p-1}}.
 $$
 In the case where $p=p_\theta$
 there exists $C_2=C_2(N,\theta)>0$ such that
  $$
  \sup_{x\in{{  \R}}^N}\mu(B(x,\sigma))
  \le C_2\left[\log\biggr(e+\frac{T^{\frac{1}{\theta}}}{\sigma}\biggr)\right]^{-\frac{N}{\theta}},
  \quad
  0<\sigma\le T^{\frac{1}{\theta}}.
  $$
 (See \cite{BP} for the  Laplacian case $\theta=2$.)
\end{itemize}
Condition~(i) implies the following non-existence result. 
  \begin{itemize}
  \item[(ii)] Let $p\ge p_\theta$. 
  There exists  
  $\gamma=\gamma  (N,\theta,p)>0$
   such that if $\mu\in{\mathcal L}_0$ satisfies  
  $$
  \begin{array}{ll}
  \mu(x)\ge \gamma|x|^{-N}\displaystyle{\biggr[\log\left(e+\frac{1}{|x|}\right)\biggr]^{-\frac{N}{\theta}-1}}
  \quad & \mbox{if}\quad \displaystyle{p=p_\theta},\vspace{7pt}\\
  \mu(x)\ge \gamma|x|^{-\frac{\theta}{p-1}}\quad & \mbox{if}\quad \displaystyle{p>p_\theta},
  \end{array}
  $$
  for almost all (a.a.) $x$ in a neighborhood of the origin, 
  then problem~\eqref{eq:1.1} possesses no local-in-time nonnegative solutions.
  \end{itemize}
Regarding sufficiency, in \cite{HI01} they obtained results (iii) and (iv) below.
  \begin{itemize}
  \item[(iii)] Let $\mu\in{\mathcal M}$ and $1<p<p_\theta$. 
  There exists 
  $c=c(N,\theta,p)>0$ such that if
 $$
 \sup_{x\in{{  \R}}^N}\mu(B(x,T^{\frac{1}{\theta}}))\le c T^{\frac{N}{\theta}-\frac{1}{p-1}}
 $$
 for some $T>0$, then problem~\eqref{eq:1.1} possesses a nonnegative solution in $Q_T$.
\item[(iv)] Let $\mu\in{\mathcal L}_0$ and $p\ge p_\theta$.
There exists $\varepsilon =\varepsilon (N,\theta,p) >0$
 such that if 
$$
\begin{array}{ll}
0\le\mu(x)\le \varepsilon |x|^{-N}\displaystyle{\biggr[\log\left(e+\frac{1}{|x|}\right)\biggr]^{-\frac{N}{\theta}-1}}+K
\quad & \mbox{if}\quad \displaystyle{p=p_\theta},\vspace{7pt}\\
0\le\mu(x)\le \varepsilon |x|^{-\frac{\theta}{p-1}}+K\quad & \mbox{if}\quad \displaystyle{p>p_\theta},
\end{array}
$$
for a.a.~$x\in {\mathbb R}^N$ for some $K>0$, then problem~\eqref{eq:1.1} possesses a local-in-time nonnegative solution.
\end{itemize}

For $\mu\in{\mathcal L}_0$ 
the results in (ii) and (iv) demonstrate that the `strength' of the singularity  at the origin of the functions
$$
\mu_c(x)=
\left\{
\begin{array}{ll}
 |x|^{-\frac{\theta}{p-1}} & \quad\mbox{if}\quad p>p_\theta,\vspace{3pt}\\
 |x|^{-N}|\log|x||^{-\frac{N}{\theta}-1} & \quad\mbox{if}\quad p=p_\theta,
\end{array}
\right.
$$
is the critical threshold 
for the local solvability of problem~\eqref{eq:1.1}. 
We term such a singularity  in the initial  data an {\em optimal singularity}
for the solvability for problem~\eqref{eq:1.1}.
Of course, by translation invariance the singularity could be located at any point of ${{  \R}}^N$. 

Subsequently, the results of \cite{HI01} were extended
to some related parabolic problems with a~power  law nonlinearity (see \cites{HI02, HIT, HS, HT, IKO, IKO2}).
However, one cannot apply the arguments in these papers to
problem~(P) with a general nonlinearity $F$ since
they depend heavily  upon the homogeneous structure of the power  law nonlinearity. 
 \subsection{The main result}
\label{subsection:1.2}
In this paper we improve the arguments in \cite{HI01} to obtain necessary conditions
and sufficient conditions for the existence of local-in-time solutions of problem~(P)
 for a significantly larger class of nonlinearities  $F$ 
 and determine the optimal singularities of the initial data for the solvability of problem~(P). 
 
Let $f_1$ and $f_2$ be real-valued functions defined in an interval $[L,\infty)$, where $L\in{{  \R}}$. We write  $f_1(t)\preceq f_2(t)$ as $t\to \infty$ if there exists $C>0$ such that 
 $f_1(t)\le Cf_2(t)$  for all large enough $t\in[L,\infty)$. We define $\succeq $ in the obvious way, namely $f_2(t)\succeq f_1(t)$ as $t\to \infty$ if and only if
  $f_1(t)\preceq f_2(t)$ as $t\to \infty$. 
 We  write $f_1(t)\asymp f_2(t)$ as $t\to \infty$ whenever $f_1(t)\preceq f_2(t)$ and $f_1(t)\succeq f_2(t)$ as $t\to \infty$,
  i.e. there exists $C>0$ such that $C^{-1}f_2(t)\le f_1(t)\le Cf_2(t)$ for large enough $t\in[L,\infty)$.
  
 We consider  nonlinearities which are asymptotic to  the prototypical example~(F), in this sense:
\begin{itemize}
  \item[{\rm (F1)}]
  $F$ is locally Lipschitz continuous in $[0,\infty)$;
  \item[{\rm (F2)}]
  $F(\tau)\asymp \tau^p[\log\tau]^q$ as $\tau\to\infty$ for some $p>1$ and $q\in{{  \R}}$.
\end{itemize}
\begin{theorem}
\label{Theorem:1.1}
Assume conditions~{\rm (F1)} and {\em (F2)}.
\begin{itemize}
  \item[{\rm (i)}]
  Let $\mu\in{\mathcal M}$ and either
  $$
  \mbox{{\rm (i)}}\quad\mbox{$1<p<p_\theta$}\qquad\mbox{or}\qquad
  \mbox{{\rm (ii)}}\quad\mbox{$p=p_\theta$ and $q<-1$}.
  $$
  Problem~{\rm (P)} possesses a local-in-time solution if and only if
  $\displaystyle{\sup_{z\in{{  \R}}^N}\mu(B(z,1))<\infty}$.
   \item[{\rm (ii)}] Suppose $\mu\in{\mathcal L}_0$. 
   \begin{itemize}
   \item[{\rm (1)}]
  Let $p=p_\theta$ and $q=-1$.
  There exists $\gamma_1>0$ such that if
  \begin{equation}
  \label{eq:1.3}
  \mu(x)\ge\gamma_1|x|^{-N}|\log |x||^{-1}[\log|\log|x||]^{-\frac{N}{\theta}-1}
  \end{equation}
  in a neighborhood of $x=0$, then problem~{\rm (P)} possesses no local-in-time solutions.
  On the other hand, 
  for any $R\in(0,1)$, 
  there exists $\varepsilon_1>0$ such that if
  \begin{equation}
  \label{eq:1.4}
  0\le\mu(x)\le\varepsilon_1|x|^{-N}|\log |x||^{-1}[\log|\log|x||]^{-\frac{N}{\theta}-1}\chi_{B(0,R)}(x)+K_1,
  \quad x\in{{  \R}}^N,
  \end{equation}
  for some $K_1>0$, then problem~{\rm (P)} possesses a local-in-time solution.
  \item[{\rm (2)}] Let $p=p_\theta$ and $q>-1$.
  There exists $\gamma_2>0$ such that if
  \begin{equation}
  \label{eq:1.5}
  \mu(x)\ge\gamma_2|x|^{-N}|\log |x||^{-\frac{N(q+1)}{\theta}-1}
  \end{equation}
  in a neighborhood of $x=0$, then problem~{\rm (P)} possesses no local-in-time solutions.
  On the other hand, for any $R\in(0,1)$, 
  there exists $\varepsilon_2>0$ such that if
  \begin{equation}
  \label{eq:1.6}
  0\le\mu(x)\le\varepsilon_2|x|^{-N}|\log |x||^{-\frac{N(q+1)}{\theta}-1}\chi_{B(0,R)}(x)+K_2,
  \quad x\in{{  \R}}^N,
  \end{equation}
  for some $K_2>0$,
  then problem~{\rm (P)} possesses a local-in-time solution.
  \item[{\rm (3)}] Let $p>p_\theta$.
  There exists $\gamma_3>0$ such that if
  \begin{equation}
  \label{eq:1.7}
  \mu(x)\ge\gamma_3|x|^{-\frac{\theta}{p-1}}|\log |x||^{-\frac{q}{p-1}}
  \end{equation}
  in a neighborhood of $x=0$, then problem~{\rm (P)} possesses no local-in-time solutions.
  On the other hand, for any $R\in(0,1)$, 
  there exists $\varepsilon_3>0$ such that if
  \begin{equation}
  \label{eq:1.8}
  0\le\mu(x)\le\varepsilon_3|x|^{-\frac{\theta}{p-1}}|\log |x||^{-\frac{q}{p-1}}\chi_{B(0,R)}(x)+K_3,
    \quad x\in{{  \R}}^N,
  \end{equation}
  for some $K_3>0$,
  then problem~{\rm (P)} possesses a local-in-time solution.
  \end{itemize}
\end{itemize}
\end{theorem}
While Theorem~\ref{Theorem:1.1} provides sharp results on the identification of optimal singularities for the solvability of problem~(P), 
we  point out that we have obtained several other interesting and powerful results in this paper 
regarding necessary conditions and sufficient conditions for existence under very general conditions on $F$. 
We mention, in particular, Theorem~\ref{Theorem:3.1},  Theorem~\ref{Theorem:4.1}, Theorem~\ref{Theorem:4.2}, and 
Theorem~\ref{Theorem:4.3}. 

Subject to mild assumptions on $F$ (essentially that of majorizing a convex function with suitable monotonicity properties), 
we follow the strategy in \cite{HI01} and obtain necessary conditions for 
the existence in Theorem~\ref{Theorem:3.1}.
However, the iteration step in \cite{HI01} to obtain the estimate for the optimal singularity
relies on the homogeneity of the pure power law nonlinearity considered there. 
For the class of nonlinearities satisfying (F1)--(F2), we combine the arguments in \cite{HI01} with the method introduced in \cite{LS02},
to obtain a sharper  necessary condition in Corollary~\ref{Corollary:3.1}. 
Conversely, in order to derive sharp sufficient conditions we require delicate arguments for $F$ satisfying (F1)--(F2).
Indeed, the arguments are separated into three cases:
(i) $1<p<p_\theta$ (see Theorem~\ref{Theorem:4.1}), $p>p_\theta$ (see Theorem~\ref{Theorem:4.3}), 
and (iii) $p=p_\theta$ (see Theorem~\ref{Theorem:4.2}). 
The arguments in case~(i) are somewhat standard but the other cases involve certain intricacies,
in particular, for the critical case $p=p_\theta$.
\vspace{3pt}

The rest of this paper is organized as follows.
In Section~\ref{section:2} we recall some properties of the fundamental solution $\Gamma_\theta$
and prove some preliminary lemmas.
In Section~\ref{section:3} we obtain necessary conditions for the existence of local-in-time solutions of problem~(P).
In Section~\ref{section:4} we prove several theorems on sufficient conditions
for the existence of local-in-time solutions of problem~(P).
In Section~\ref{subsection:4.4} 
we also provide a necessary and sufficient condition on the nonlinearity $F$ for which problem~(P) is solvable 
for the case  of initial data a Dirac measure   (Corollary~\ref{Corollary:4.4}). 
Finally, in Section~\ref{section:5}
we complete the proof of our main result, Theorem~\ref{Theorem:1.1}, 
and outline some analogous results
for nonlinearities which are asymptotic to further log-refinements of the cases above (see Remark~\ref{Remark:5.1}).
\section{Preliminaries}
\label{section:2}
In this section we prove some important technical lemmas, 
modifying the arguments in \cite{HI01} for the more general nonlinearities considered here.
We  make precise our notion of solution used throughout this paper, which implicitly considers nonnegative functions only. 
The use of the word `solvability' for problem~(P) is always used with respect to this solution concept.
In all that follows we will use $C$ to denote generic positive constants which depend only on $N$, $\theta$, and~$F$ and point out that $C$  
 may take different values  within a calculation.
We begin by recalling some properties of the kernel for the fractional Laplacian.  

Let $\Gamma_\theta=\Gamma_\theta(x,t)$ be the fundamental solution of
$$
\partial_t u+(-\Delta)^{\frac{\theta}{2}}u=0\quad\mbox{in}\quad{{  \R}}^N\times(0,\infty).
$$
The function $\Gamma_\theta$ satisfies
\begin{equation}
\label{eq:2.1}
\begin{split}
 & \Gamma_\theta(x,t)=(4\pi t)^{-\frac{N}{2}}\exp\left(-\frac{|x|^2}{4t}\right)\quad\mbox{if}\quad \theta=2,\\
 & C^{-1}t^{-\frac{N}{\theta}}\left(1+t^{-\frac{1}{\theta}}|x|\right)^{-N-\theta}\le\Gamma_\theta(x,t)
\le Ct^{-\frac{N}{\theta}}\left(1+t^{-\frac{1}{\theta}}|x|\right)^{-N-\theta}\quad\mbox{if}\quad 0<\theta<2,
\end{split}
\end{equation}
for all $x\in{{  \R}}^N$ and $t>0$ and  has the following properties:
\begin{eqnarray}
\nonumber
 & \bullet & \mbox{$\Gamma_\theta$ is positive and smooth in ${{  \R}}^N\times(0,\infty)$},\\
\label{eq:2.2}
 & \bullet & \Gamma_\theta(x,t)=t^{-\frac{N}{\theta}}\Gamma_\theta\left(t^{-\frac{1}{\theta}}x,1\right),
 \quad \int_{{\mathbb R}^N}\Gamma_\theta(x,t)\,dx=1,\\
\label{eq:2.3}
 & \bullet &  \mbox{$\Gamma_\theta(\cdot,1)$ is radially symmetric and $\Gamma_\theta(x,1)\le \Gamma_\theta(y,1)$ if $|x|\ge |y|$},\\
\label{eq:2.4}
 & \bullet & \Gamma_\theta(x,t)=\int_{  \R^N}\Gamma_\theta(x-y,t-s)\Gamma_\theta(y,s)\,dy,
\end{eqnarray}
for all $x,y\in{{  \R}}^N$ and $0<s<t$ 
(see for example \cites{BJ, BraK, S}). 
Furthermore, we have the following smoothing estimate
for  the semigroup associated with $\Gamma_\theta$
(see \cite{HI01}*{Lemma~2.1}).
\begin{lemma}
\label{Lemma:2.1}
For any $\mu\in{\mathcal M}$, set
$$
[S(t)\mu](x):= \int_{{{  \R}}^N} \Gamma_\theta(x-y,t)\,d\mu(y), \quad x\in{{  \R}}^N,\,\,\, t>0.
$$
Then there exists $C=C(N,\theta)>0$ such that
$$
\|S(t)\mu\|_{L^\infty({{  \R}}^N)} \le C t^{-\frac{N}{\theta}} \sup_{x\in{{  \R}}^N} \mu(B(x,t^\frac{1}{\theta})),
\quad t>0.
$$
\end{lemma}
\begin{remark}
\label{Remark:2.1}
{\rm (i)} $S(t)\mu$ is possibly infinite everywhere in ${ \R}^N$; {\rm (ii)} if $\mu\in{\mathcal M}$ is such that 
$$
\sup_{x\in{ \R}^N}\mu(B(x,r ))<\infty
$$
for some $r>0$, then for any $R\ge r$ there exists $C\ge 1$ such that 
$$
\sup_{x\in{ \R}^N}\mu(B(x,R))\le C\sup_{x\in{ \R}^N}\mu(B(x,r ))<\infty.
$$
See for example \cite{IS01}*{Lemma~2.1} or  \cite{FI01}*{Lemma~2.4}.
\end{remark}

We now make precise our  solution concepts for problem~(P).
\begin{definition}
\label{Definition:2.1}
Let $T>0$ and $u$ be a nonnegative, measurable, finite almost everywhere function in $Q_T$. 
Let $F$ be a nonnegative and continuous function in $[0,\infty)$. 
\begin{itemize}
  \item[{\rm (i)}]
  We say that $u$ satisfies
  \begin{equation}
  \label{eq:2.5}
  \partial_t u+(-\Delta)^{\frac{\theta}{2}}u=F(u)
  \end{equation}
  in $Q_T$ if, for a.a.~$\tau\in(0,T)$, $u$ satisfies
  $$
  u(x,t)=\int_{{{  \R}}^N}\Gamma_\theta(x-y,t-\tau)u(y,\tau)\,dy+\int_\tau^t\int_{{{  \R}}^N}\Gamma_\theta(x-y,t-s)F(u(y,s))\,dy\,ds
  $$
  for a.a.~$(x,t)\in{{  \R}}^N\times (\tau,T)$. 
  \item[{\rm (ii)}]
  Let $\mu\in{\mathcal M}$. 
  We say that $u$ is a solution of problem~{\rm (P)}
  in $Q_T$ if $u$ satisfies
  \begin{equation}
  \label{eq:2.6}
  u(x,t)=\int_{{{  \R}}^N}\Gamma_\theta(x-y,t)\,d\mu(y)+\int_0^t\int_{{{  \R}}^N}\Gamma_\theta(x-y,t-s)F(u(y,s))\,dy\,ds
  \end{equation}
  for a.a.~$(x,t)\in Q_T$. 
  If $u$ satisfies \eqref{eq:2.6} with $``="$ replaced by $``\ge"$,
  then $u$ is said to be a supersolution of problem~{\rm (P)} in $Q_T$.
\end{itemize}
\end{definition}
Next we recall a lemma on the existence of solutions of problem~(P) in the presence of a~supersolution 
(see \cite{HI01}*{Lemma~2.2}). 
\begin{lemma}
\label{Lemma:2.2}
Let 
$F$ be an increasing, nonnegative continuous function in $[0,\infty)$.
Let $\mu\in{\mathcal M}$ and $0<T\le\infty$. 
If there exists a supersolution~$v$ of problem~{\rm (P)} in $Q_T$, 
then there exists a solution~$u$ of problem~{\rm (P)} in $Q_T$ such that
$0\le u(x,t)\le v(x,t)$ in $Q_T$.
\end{lemma}
Combining Lemma~\ref{Lemma:2.2} and parabolic regularity theory,
we have:
\begin{lemma}
\label{Lemma:2.3}
Let $\mu\in{\mathcal M}$ be such that
$\displaystyle{\sup_{z\in{{  \R}}^N}}\mu(B(z,1))<\infty$. Suppose
\begin{itemize}
  \item[{\rm (i)}]
   $F_1$ is nonnegative and locally Lipschitz continuous in $[0,\infty)$; 
  \item[{\rm (ii)}]
   $F_2$ is an increasing and  continuous function in $[0,\infty)$ such that
  $F_1(\tau)\le F_2(\tau)$ for all $\tau\in[0,\infty)$.
  \end{itemize}
  If there exists a supersolution $v$ of {\rm (P)}  in $Q_T$ with $F$ replaced by $F_2$
   such that for all $\tau\in(0,T)$
  \begin{equation}
  \label{eq:2.7}
  \sup_{\tau<t<T}\|v(t)\|_{L^\infty({{  \R}}^N)}<\infty ,
  \end{equation}
then there exists a solution~$u$ of {\rm (P)} in $Q_T$ with $F$ replaced by $F_1$, 
with $u$ satisfying $0\le u(x,t)\le v(x,t)$ in $Q_T$. 
\end{lemma}
{\bf Proof.}
For any $m$, $n\in\N$ set
\begin{eqnarray*}
F_{1,m}(\tau)&:=& \min\{F_1(\tau),m\}\quad\mbox{ for $\tau\ge 0$},\\
\mu_n(x)& :=& \int_{{{  \R}}^N}\Gamma_\theta(x-y,2n^{-1})\,d\mu(y)\quad\mbox{ for $x\in{{  \R}}^N$}.
\end{eqnarray*}
It follows from Lemma~\ref{Lemma:2.1} that $S(n^{-1})\mu\in L^\infty({{  \R}}^N)$.
Also, since $\mu_n=S(n^{-1})S(n^{-1})\mu$, we have that $\mu_n\in BC({{  \R}}^N)$.
For each $m$, $n\in \N$ define the sequence $\{u_{m,n,k}\}_{k=0}^\infty$ by
\begin{equation*}
\begin{split}
 & u_{m,n,0}(x,t):=\int_{{{  \R}}^N}\Gamma_\theta(x-y,t)\mu_n(y)\,dy,\\
 & u_{m,n,k+1}(x,t):=\int_{{{  \R}}^N}\Gamma_\theta(x-y,t)\mu_n(y)\,dy
 +\int_0^t\int_{{{  \R}}^N}\Gamma_\theta(x-y,t-s)F_{1,m}(u_{m,n,k}(y,s))\,dy\,ds.
\end{split}
\end{equation*}
By \eqref{eq:2.4} and Definition~\ref{Definition:2.1}~(ii) we have 
\begin{align*}
u_{m,n,0}(x,t) & =\int_{{{\R}}^N}\Gamma_\theta(x-y,t)\left(\int_{{{\R}}^N}\Gamma_\theta(y-z,2n^{-1})\,d\mu(z)\right)\,dy\\
 & =\int_{{{\R}}^N}\Gamma_\theta(x-z,t+2n^{-1})\,d\mu(z)\\
 & \le v(x,t+2n^{-1})
\end{align*}
for $x\in{{  \R}}^N$ and $t\in[0,T-2n^{-1})$. 
Since $F_1(\tau)\le F_2(\tau)$ for $\tau\in[0,\infty)$,
by induction we obtain
\begin{equation}
\label{eq:2.8}
0\le u_{m,n,k}(x,t)\le v(x,t+2n^{-1})
\end{equation}
for all $x\in{{  \R}}^N$, $t\in[0,T-2n^{-1})$, and $k\ge 0$.
Here we used the assumption that $F_2$ is increasing.
Since $F_{1,m}$ is globally Lipschitz in $[0,\infty)$,
we may apply the standard theory of evolution equations to see that the pointwise limit
$$
u_{m,n}(x,t):=\lim_{k\to\infty}u_{m,n,k}(x,t)
$$
exists in ${{  \R}}^N\times[0,\infty)$ and  satisfies
\begin{equation}
\label{eq:2.9}
u_{m,n}(x,t)=\int_{{{  \R}}^N}\Gamma_\theta(x-y,t)\mu_n(y)\,dy+\int_0^t\int_{{{  \R}}^N}\Gamma_\theta(x-y,t-s)F_{1,m}(u_{m,n}(y,s))\,dy\,ds
\end{equation}
for all $x\in{{  \R}}^N$ and $t>0$.
Furthermore, by \eqref{eq:2.8} we see that 
\begin{equation}
\label{eq:2.10}
0\le u_{m,n}(x,t)\le v(x,t+2n^{-1})
\end{equation}
for all $x\in{{  \R}}^N$ and $t\in[0,T-2n^{-1})$.  
Then, by \eqref{eq:2.7}, for any $\tau\in(0,T-2n^{-1})$ we have 
$$
\sup_{\tau<t<T-2n^{-1}}\|u_{m,n}(t)\|_{L^\infty({{  \R}}^N)}\le\sup_{\tau<t<T}\|v(t)\|_{L^\infty({{  \R}}^N)}<\infty.
$$
Applying the standard
parabolic regularity theory to integral equation~\eqref{eq:2.9},
we find $\alpha\in(0,1)$ such that
\begin{equation}
\label{eq:2.11}
\sup_n \|u_{m,n}\|_{C^{\alpha;\alpha/2}(K)}<\infty
\end{equation}
for any compact set $K\subset Q_T$.
By the Ascoli-Arzel\'a theorem and the  diagonal argument
we obtain a subsequence $\{u_{m,n'}\}$ of $\{u_{m,n}\}$ and a function $u_m\in C(Q_T)$ 
such that 
\begin{equation}
\label{eq:2.12}
\lim_{n'\to\infty}u_{m,n'}(x,t)=u_m(x,t)\quad\mbox{in}\quad Q_T.
\end{equation}
Since $F_{1,m}$ is bounded and continuous in $(0,\infty)$, 
by \eqref{eq:2.9} we have
\begin{equation}
\label{eq:2.13}
u_m(x,t)=\int_{{{  \R}}^N}\Gamma_\theta(x-y,t)\mu(y)\,dy+\int_0^t\int_{{{  \R}}^N}\Gamma_\theta(x-y,t-s)F_{1,m}(u_m(y,s))\,dy\,ds
\end{equation}
in $Q_T$.
Furthermore, by \eqref{eq:2.10} and \eqref{eq:2.12} we see that 
\begin{equation}
\label{eq:2.14}
0\le u_m(x,t)\le v(x,t)
\end{equation}
for all $x\in{{  \R}}^N$ and $t\in[0,T)$.  

Similarly to \eqref{eq:2.11}, 
using \eqref{eq:2.14}, instead of \eqref{eq:2.10},
we have 
$$
\sup_n \|u_m\|_{C^{\alpha;\alpha/2}(K)}<\infty
$$
for any compact set $K\subset Q_T$.
By the Ascoli-Arzel\'a theorem and the diagonal argument
we obtain a subsequence $\{u_{m'}\}$ of $\{u_m\}$ and a function $u\in C(Q_T)$ 
such that 
\begin{equation}
\label{eq:2.15}
\lim_{m'\to\infty}u_{m'}(x,t)=u(x,t)\quad\mbox{in}\quad Q_T.
\end{equation}
Since $F_{1,m}(\tau)\le F_1(\tau)\le F_2(\tau)$ for $\tau\in(0,\infty)$, 
by \eqref{eq:2.14} we see that 
\begin{align*}
 & \sup_{m'}\int_0^t\int_{{{  \R}}^N}
\Gamma_\theta(x-y,t-s)F_{1,m'}(u_{m'}(y,s))\,dy\,ds\\
 & \le\int_0^t\int_{{{  \R}}^N}
\Gamma_\theta(x-y,t-s)F_2(v(y,s))\,dy\,ds
\le v(x,t)<\infty
\end{align*}
for a.a.~$(x,t)\in Q_T$. 
Then, by \eqref{eq:2.13} and \eqref{eq:2.15}
we apply Lebesgue's dominated convergence theorem to see that $u$ is a solution of problem~(P) in $Q_T$ with $F$ replaced by $F_1$ 
and $0\le u(x,t)\le v(x,t)$ in $Q_T$. 
Thus Lemma~\ref{Lemma:2.3} follows.
\qed
\vspace{5pt}

Next we provide two lemmas
on the relationship between the initial measure and the initial trace for problem~(P).
\begin{lemma}
\label{Lemma:2.4}
Let $F$ be a nonnegative continuous function in $[0,\infty)$.
\begin{itemize}
  \item[{\rm (i)}]
  Let $u$ satisfy \eqref{eq:2.5} in $ Q_T$ for some $T>0$. Then
  $$
  \underset{0<t<T-\varepsilon}{\mbox{{\rm ess sup}}}\,\int_{B(0,R)}u(y,t)\,dy<\infty
  $$
  for all $R>0$ and $0<\varepsilon<T$.
  Furthermore,
  there exists a unique $\nu\in{\mathcal M}$ as an initial trace of the solution~$u$; that is,
  \begin{equation*}
  \underset{t\to +0}{\mbox{{\rm ess lim}}}
  \int_{{{  \R}}^N}u(y,t)\eta(y)\,dy=\int_{{{  \R}}^N}\eta(y)\,d\nu(y)
  \end{equation*}
  for all $\eta\in C_0({{  \R}}^N)$.
  \item[{\rm (ii)}]
  Let $u$ be a solution of problem~{\rm (P)} in $Q_T$ for some $T>0$.
  Then assertion~{\rm (i)} holds with $\nu$ replaced by $\mu$.
\end{itemize}
\end{lemma}
The proof of Lemma~\ref{Lemma:2.4}
is the same as in 
\cite{HI01}*{Lemma~2.3}.
Furthermore, by assertion~(i)
we can apply the same argument as in the proof of \cite{HI01}*{Theorem~1.2} to obtain the following lemma.
\begin{lemma}
\label{Lemma:2.5}
Let $F$ be a  nonnegative continuous function in $[0,\infty)$ and $T>0$. 
Let $u$ satisfy \eqref{eq:2.5} in $Q_T$. 
Let $\mu\in{\mathcal M}$ be the unique  initial trace of $u$  guaranteed by Lemma~{\rm\ref{Lemma:2.4}}.
If $\sup_{z\in{{  \R}}^N}\mu(B(z,1))<\infty$ then $u$ is a solution of problem~{\rm (P)} in $Q_T$. 
\end{lemma}

In the rest of this section 
we prepare preliminary lemmas. 
\begin{lemma}
\label{Lemma:2.6}
Let $a>0$ and $b$, $c\in{\mathbb R}$. 
Set 
$$
\varphi(\tau):=\tau^a(\log \tau)^b(\log\log \tau)^c,\quad \tau\in(e,\infty). 
$$
Then there exists $L\in(e,\infty)$ such that $\varphi'>0$ in $(L,\infty)$ 
and  the inverse function $\varphi^{-1}:(\varphi(L),\infty)\to (L,\infty)$ exists. 
Furthermore,
\begin{equation}
\label{eq:2.16}
\varphi^{-1}(\tau)\asymp\tau^{\frac{1}{a}}(\log\tau)^{-\frac{b}{a}}(\log\log\tau)^{-\frac{c}{a}}
\end{equation}
as $\tau\to\infty$. 
\end{lemma}
{\bf Proof.}
Since $a>0$, we can find $L\in(e,\infty)$ such that 
$$
\varphi'(\tau)=
\tau^{a-1}(\log\tau)^b(\log\log \tau)^c
\left[a+b(\log\tau)^{-1}+c(\log\tau)^{-1}(\log\log \tau)^{-1}\right]>0
$$
 for all $\tau\in(L,\infty)$. 
 Since $\varphi (\tau)\to\infty$ as $\tau\to\infty$, it follows that $\varphi^{-1}:(\varphi(L),\infty)\to (L,\infty)$  exists and 
 satisfies $\varphi^{-1}(\tau)\to\infty$ as $\tau\to\infty$.
Now, 
\begin{eqnarray}
\log \tau & = &  \log \varphi(\varphi^{-1}(\tau))=a\log\varphi^{-1}(\tau)+b\log\log\varphi^{-1}(\tau)+c\log\log\log\varphi^{-1}(\tau)\label{eq:2.17}\\
 & = & a\log\varphi^{-1}(\tau)(1+o(1))\notag
\end{eqnarray}
as $\tau\to\infty$, so that
$$
\log\varphi^{-1}(\tau)=\frac{1}{a}(\log \tau)(1+o(1))
$$
as $\tau\to\infty$.
Then, by \eqref{eq:2.17} we have
\begin{equation*}
\begin{split}
a\log\varphi^{-1}(\tau) & =\log\tau-b\log\log\varphi^{-1}(\tau)-c\log\log\log\varphi^{-1}(\tau)\\
 & =\log \tau -b\log\left(\frac{1}{a}(\log \tau)(1+o(1))\right)-c\log\log\left(\frac{1}{a}(\log \tau)(1+o(1))\right)
\end{split}
\end{equation*}
as $\tau\to\infty$. Hence, 
\begin{equation*}
\log\varphi^{-1}(\tau)   =\log\left[\left(\tau  (\log\tau)^{-b}(\log\log\tau)^{-c}\right)^{\frac{1}{a}}\right]-\frac{b}{a}\log\left(\frac{1}{a}(1+o(1))\right)+o(1)
\end{equation*}
as $\tau\to\infty$, from which  \eqref{eq:2.16} follows  and completes the proof of Lemma~\ref{Lemma:2.6}. 
\qed
\begin{lemma}
\label{Lemma:2.7}
Let $a>0$, $b\ge 0$, and $c\in{  \R}$. 
\vspace{3pt}
\newline
{\rm (i)} There exists $C_1>0$ such that 
$$
\int_A^B\tau^{a-b-1}(\log \tau)^c\,d\tau
\ge C_1 A^aB^{-b}(\log A)^c\log\frac{B}{A}
$$
for all $A$, $B\in[2,\infty)$ with $A\le B$. 
\vspace{3pt}
\newline
{\rm (ii)} There exists $C_2\in[1,\infty)$ such that 
$$
C_2^{-1}\tau^a[\log(e+\tau)]^c\le
\int_{\tau/2}^\tau s^{a-1}[\log(e+s)]^c\,ds
\le \int_0^\tau s^{a-1}[\log(e+s)]^c\,ds
\le C_2\tau^a[\log(e+\tau)]^c
$$
for all $\tau\in[0,\infty)$. 
\end{lemma}
{\bf Proof.}
We first prove assertion~(i). 
Thanks to $a>0$, 
by Lemma~\ref{Lemma:2.6}
we find $R_1\in[2,\infty)$ such that 
the function $(R_1,\infty)\ni\tau\mapsto \tau^a(\log \tau)^c$ is increasing. 
Then we have
\begin{equation*} 
\left\{
\begin{array}{ll}
\tau^a(\log \tau)^c\ge A^a(\log A)^c & \quad\mbox{if}\quad R_1\le A\le\tau,\vspace{3pt}\\
\tau^a(\log \tau)^c \ge R_1^a(\log R_1)^c\ge CA^a(\log A)^c & \quad\mbox{if}\quad A\le R_1\le\tau,\vspace{3pt}\\
\tau^a(\log \tau)^c \ge CA^a(\log A)^c & \quad\mbox{if}\quad A\le\tau<R_1,\qquad
\end{array}
\right.
\end{equation*}
for all $A$, $B\in[2,\infty)$ with $A\le B$. 
We notice that
$$
\inf_{A\in[2,R_1]}\frac{R_1^a(\log R_1)^c}{A^a(\log A)^c}>0,
\qquad
\inf_{A\in[2,R_1],\,\tau\in[A,R_1]}\frac{\tau^a(\log \tau)^c}{A^a(\log A)^c}>0. 
$$
We observe that
\begin{equation*}
\begin{split}
\int_A^B
\tau^{a-b-1}(\log \tau)^c\,d\tau
 & \ge B^{-b}
\int_A^B\tau^{a-1}(\log\tau)^c\,d\tau\\
 & \ge CB^{-b}A^a(\log A)^c\int_A^B\tau^{-1}\,d\tau
 =CA^aB^{-b}(\log A)^c\log\frac{B}{A}
\end{split}
\end{equation*}
for all $A$, $B\in[2,\infty)$ with $A\le B$. 
Then assertion~(i) follows. 

Next we prove assertion~(ii). 
Let $\varepsilon\in(0,a)$. 
By Lemma~\ref{Lemma:2.6} 
we find $R_2>0$ such that 
the function $(R_2,\infty)\ni\tau\mapsto (e+\tau)^\varepsilon(\log(e+\tau))^c$ is increasing. 
Then we have 
\begin{align*}
\int_0^\tau s^{a-1}[\log(e+s)]^c\,ds
 & \le C+\int_{R_2}^\tau s^{a-1}(e+s)^{-\varepsilon}(e+s)^\varepsilon[\log(e+s)]^c\,ds\\
 & \le C+(e+\tau)^\varepsilon[\log(e+\tau)]^c
\int_{R_2}^\tau s^{a-1-\varepsilon}\,ds\\
& \le C+C\tau^{a-\varepsilon}(e+\tau)^\varepsilon[\log(e+\tau)]^c
\le C\tau^a[\log(e+\tau)]^c
\end{align*}
for all $\tau\in[R_2,\infty)$. 
On the other hand, 
$$
\int_0^\tau s^{a-1}[\log(e+s)]^c\,ds
\le C\int_0^\tau s^{a-1}\,ds\le C\tau^a\le C\tau^a[\log(e+\tau)]^c
$$
for all $\tau\in(0,R_2)$. 
These imply that 
\begin{equation}
\label{eq:2.18}
\int_0^\tau s^{a-1}[\log(e+s)]^c\,ds\le C\tau^\varepsilon[\log(e+\tau)]^c\int_0^\tau s^{a-\varepsilon-1}\,ds
\le C\tau^a[\log(e+\tau)]^c
\end{equation}
for all $\tau\in[0,\infty)$. 
On the other hand, 
since
$$
\inf_{\tau\in(0,\infty)}\frac{\log(e+\tau/2)}{\log(e+\tau)}>0, 
$$
we have
$$
C^{-1}\log(e+\tau)\le \log(e+\tau/2)\le \inf_{\xi\in(\tau/2,\tau)}\log(e+\xi)\le \sup_{\xi\in(\tau/2,\tau)}\log(e+\xi)\le\log(e+\tau)
$$
for $\tau>0$. This yields 
\begin{equation}
\label{eq:2.19}
\int_{\tau/2}^\tau s^{a-1}[\log(e+s)]^c\,ds\ge C[\log(e+\tau)]^c\int_{\tau/2}^\tau s^{a-1}\,ds
\ge C\tau^a[\log(e+\tau)]^c
\end{equation}
for all $\tau\in[0,\infty)$. 
By \eqref{eq:2.18} and \eqref{eq:2.19} we have assertion~(ii). The proof is complete.
\qed
\begin{lemma}
\label{Lemma:2.8}
Let $p>1$, $d\in [1,p)$, $q\in{\mathbb R}$, and $R\ge 0$. 
Define a function~$f$ in $[0,\infty)$ by 
$$
f(\tau):=
\left\{
\begin{array}{ll}
0 & \quad\mbox{for}\quad \tau\in[0,R],\vspace{5pt}\\
\displaystyle{\tau^d\int_R^\tau s^{-d}\left(\int_R^s \xi^{p-2}[\log(e+\xi)]^q\,d\xi\right)\,ds} & \quad\mbox{for}\quad \tau\in(R,\infty).
\end{array}
\right.
$$
Then 
\begin{itemize}
  \item[{\rm (i)}] the function $(0,\infty)\ni\tau\mapsto \tau^{-d}f(\tau)$ is increasing;
  \item[{\rm (ii)}] $f$ is convex in $[0,\infty)$;
  \item[{\rm (iii)}] $f(\tau)\asymp \tau^p(\log \tau)^q$ as $\tau\to\infty$.
\end{itemize}
\end{lemma}
{\bf Proof.}
By the definition of $f$ we easily obtain property~(i). 
Since
$$
f'(\tau)=d\tau^{d-1}
\int_R^\tau s^{-d}\left(\int_R^s \xi^{p-2}[\log(e+\xi)]^q\,d\xi\right)\,ds
+\int_{R}^\tau \xi^{p-2}[\log(e+\xi)]^q\,d\xi
$$
for $\tau\in(R,\infty)$, we observe that $f'$ is increasing in $[0,\infty)$, so that  
property~(ii) holds. 

We prove property~(iii). 
Since $d\in (1,p)$, by Lemma~\ref{Lemma:2.7}~(ii) 
we have
\begin{equation}
\begin{split}
\label{eq:2.20}
f(\tau) & \le \displaystyle{\tau^d\int_0^\tau s^{-d}\left(\int_0^s \xi^{p-2}[\log(e+\xi)]^q\,d\xi\right)\,ds}\\
& \le C\tau^d\int_0^\tau s^{p-1-d}[\log(e+s)]^q\,ds
\le C\tau^p[\log(e+\tau)]^q
\end{split}
\end{equation}
for  all $\tau > R$ and 
\begin{equation}
\label{eq:2.21}
\begin{split}
f(\tau) & \ge \tau^d\int_{\tau/2}^\tau s^{-d}\left(\int_{s/2}^s \xi^{p-2}[\log(e+\xi)]^q\,d\xi\right)\,ds\\
 &  \ge C\tau^d\int_{\tau/2}^\tau s^{p-1-d}[\log(e+s)]^q\,ds
\ge C\tau^p[\log(e+\tau )]^q
\end{split}
\end{equation}
for all $\tau > 4R$.
By \eqref{eq:2.20} and \eqref{eq:2.21} we obtain assertion~(iii). 
Thus Lemma~\ref{Lemma:2.8} follows.
\qed
\section{Necessary Conditions for Solvability  }
\label{section:3}
In this section we establish necessary conditions for the solvability of problem~(P). 
We begin in Theorem~\ref{Theorem:3.1} by imposing only weak constraints on the nonlinearity $F$, 
before specializing to the case where $F$ satisfies (F1) and (F2) in Corollary~\ref{Corollary:3.1}. 
\begin{theorem}
\label{Theorem:3.1}
Let $F$ be a continuous function in $[0,\infty)$.
Assume that there exists a convex function $f$ in $[0,\infty)$ with the following properties:
\begin{itemize}
  \item[{\rm (f1)}]
  $F(\tau)\ge f(\tau)\ge 0$ in $[0,\infty)$;
  \item[{\rm (f2)}]
  the function $(0,\infty)\ni \tau\mapsto \tau^{-d}f(\tau)$ is increasing for some $d>1$.
\end{itemize}
Let $u$ satisfy \eqref{eq:2.5} in $Q_T$ for some $T>0$ and
let $\mu$ be the initial trace of $u$. 
Then there exists $\gamma=\gamma(N,\theta,f)\ge 1$ 
such that
\begin{equation}
\label{eq:3.1}
\int_{\gamma^{-1}T^{-\frac{N}{\theta}}m_\sigma(z)}^{\gamma^{-1}\sigma^{-N}m_\sigma(z)}
s^{-p_\theta-1}f(s)\,ds\le \gamma^{p_\theta+1} m_\sigma(z)^{-\frac{\theta}{N}}
\end{equation}
for all $z\in{{  \R}}^N$ and $\sigma\in(0,T^{\frac{1}{\theta}})$, where
$m_\sigma(z):=\mu(B(z,\sigma))$.
\end{theorem}
{\bf Proof.}
It follows from Definition~\ref{Definition:2.1}~(i) and property~(f1) that, 
for a.a.~$\tau\in(0,T)$, 
$$
\infty >u(x,t)
\ge\int_{{{  \R}}^N}\Gamma_\theta(x-y,t-\tau)u(y,\tau)\,dy
+\int_\tau^t\int_{{{  \R}}^N}\Gamma_\theta(x-y,t-s)f(u(y,s))\,dy\,ds
$$
for a.a.~$x\in{{  \R}}^N$ and a.a.~$t\in(\tau,T)$. 
This implies that
\begin{equation}
\label{eq:3.2}
\infty>u(x,2t)\ge\int_{{{  \R}}^N}\Gamma_\theta(x-y,t)u(y,t)\,dy
\end{equation}
for a.a.~$x\in{{  \R}}^N$ and a.a.~$t\in(0,T/2)$. 

Let $0<\rho<(T/2)^{\frac{1}{\theta}}$.
It follows from Definition~\ref{Definition:2.1}~(i), property~(f1), and \eqref{eq:2.4} that
\begin{equation}
\label{eq:3.3}
\begin{split}
 & \int_{{{  \R}}^N}\Gamma_\theta(z-x,t)u(x,t)\,dx\\
 & \ge\int_{{{  \R}}^N}\int_{{{  \R}}^N}\Gamma_\theta(z-x,t)\Gamma_\theta(x-y,t)\,d\mu(y)\,dx\\
 & \qquad\quad
 +\int_0^t\int_{{{  \R}}^N}\int_{{{  \R}}^N}\Gamma_\theta(z-x,t)\Gamma_\theta(x-y,t-s)f(u(y,s))\,dy\,ds\,dx\\
 & =\int_{{{  \R}}^N}\Gamma_\theta(z-y,2t)\,d\mu(y)
 +\int_0^t\int_{{{  \R}}^N}\Gamma_\theta(z-y,2t-s)f(u(y,s))\,dy\,ds
\end{split}
\end{equation}
for all~$z\in{{  \R}}^N$ and a.a.~$t\in(0,T)$. 
On the other hand, by \eqref{eq:2.1} we have
\begin{equation}
\label{eq:3.4}
\begin{split}
\int_{{{  \R}}^N}\Gamma_\theta(z-y,2t)\,d\mu(y)
 & \ge \int_{B(z,\sigma)}\Gamma_\theta(z-y,2t)\,d\mu(y)\\
 & \ge\min_{x\in B(0,\sigma)}\Gamma_\theta(x,2t) \mu(B(z,\sigma))
\ge Ct^{-\frac{N}{\theta}}\mu(B(z,\sigma))
\end{split}
\end{equation}
for all $z\in{{  \R}}^N$ and $t\ge\rho^\theta$, where $\sigma:=2^{\frac{1}{\theta}}\rho\in(0,T^{\frac{1}{\theta}})$.
Furthermore,
by \eqref{eq:2.2} and \eqref{eq:2.3} we see that
\begin{equation}
\label{eq:3.5}
\begin{split}
\Gamma_\theta(z-y,2t-s)
 & =(2t-s)^{-\frac{N}{\theta}}\Gamma_\theta\left((2t-s)^{-\frac{1}{\theta}}(z-y),1\right)\\
 & \ge\left(\frac{s}{2t}\right)^{\frac{N}{\theta}}s^{-\frac{N}{\theta}}\Gamma_\theta\left(s^{-\frac{1}{\theta}}(z-y),1\right)
=\left(\frac{s}{2t}\right)^{\frac{N}{\theta}}\Gamma_\theta(z-y,s)
\end{split}
\end{equation}
for  all $y$, $z\in{{  \R}}^N$ and $0<s<t$.
Combining \eqref{eq:3.2}, \eqref{eq:3.3}, \eqref{eq:3.4}, and \eqref{eq:3.5},
we find $C_*\ge 1$ such that
\begin{equation*}
\begin{split}
\infty>w(t):= & \int_{{{  \R}}^N}\Gamma_\theta(z-x,t)u(x,t)\,dx\\
\ge &\, C_*^{-1}t^{-\frac{N}{\theta}}\mu(B(z,\sigma))
+C_*^{-1}t^{-\frac{N}{\theta}}\int_{\rho^\theta}^t\int_{{{  \R}}^N}s^{\frac{N}{\theta}}\Gamma_\theta(z-y,s)f(u(y,s))\,dy\,ds
\end{split}
\end{equation*}
for a.a.~$z\in{{  \R}}^N$ and a.a.~$t\in(\rho^\theta,T/2)$.
Thanks to the convexity of $f$, by \eqref{eq:2.2} 
we may apply Jensen's inequality to obtain
\begin{equation}
\label{eq:3.6}
\begin{split}
\infty>w(t) & \ge C_*^{-1}t^{-\frac{N}{\theta}}\mu(B(z,\sigma))
+C_*^{-1}t^{-\frac{N}{\theta}}\int_{\rho^\theta}^t s^{\frac{N}{\theta}}
f\left(\int_{{{  \R}}^N} \Gamma_\theta(z-y,s)u(y,s)\,dy\right)\,ds\\
 & =C_*^{-1}t^{-\frac{N}{\theta}}\mu(B(z,\sigma))
 +C_*^{-1}t^{-\frac{N}{\theta}}\int_{\rho^\theta}^t s^{\frac{N}{\theta}}f(w(s))\,ds
\end{split}
\end{equation}
for a.a.~$z\in{{  \R}}^N$ and a.a.~$t\in(\rho^\theta,T/2)$.

Since $f$ is convex in $[0,\infty)$, it is Lipschitz continuous in any compact subinterval of $[0,\infty)$. We may then let
$\zeta$ denote the unique local solution of the integral equation
\begin{equation}
\label{eq:3.7}
\zeta(t)=\mu(B(z,\sigma))
+\int_{\rho^\theta}^t s^{\frac{N}{\theta}}f(C_*^{-1}s^{-\frac{N}{\theta}}\zeta(s))\,ds,\quad t\ge \rho^\theta.
\end{equation}
Hence, $\zeta$ is the unique local solution of 
\begin{equation}
\label{eq:3.8}
\zeta'(t)=t^{\frac{N}{\theta}}f(C_*^{-1}t^{-\frac{N}{\theta}}\zeta(t)),
\qquad
\zeta(\rho^\theta)=\mu(B(z,\sigma)).
\end{equation}
By \eqref{eq:3.6}, applying the standard theory for ordinary differential equations to \eqref{eq:3.7},
we see that the solution~$\zeta$ exists in $[\rho^\theta,T/2)$ and  satisfies
$$
\zeta(t)\le C_* t^{\frac{N}{\theta}}w(t)<\infty,\quad t\in[\rho^\theta,T/2).
$$
It follows from \eqref{eq:3.8} and property~(f2) that
\begin{equation*}
\begin{split}
\zeta'(t) & =t^{\frac{N}{\theta}}
[C_*^{-1}t^{-\frac{N}{\theta}}\zeta(t)]^d
[C_*^{-1}t^{-\frac{N}{\theta}}\zeta(t)]^{-d}
f\left(C_*^{-1}t^{-\frac{N}{\theta}}\zeta(t)\right)\\
 & \ge t^{\frac{N}{\theta}}
C_*^{-1}t^{-\frac{N}{\theta}}\zeta(t)]^d[C_*^{-1}t^{-\frac{N}{\theta}}\zeta(\rho^\theta)]^{-d}
f\left(C_*^{-1}t^{-\frac{N}{\theta}}\zeta(\rho^\theta)\right)\\
& \ge \zeta(\rho^\theta)^{-d}t^{\frac{N}{\theta}}\zeta(t)^d
 f\left(C_*^{-1}t^{-\frac{N}{\theta}}\zeta(\rho^\theta)\right)
\end{split}
\end{equation*}
for all $t\in[\rho^\theta,T/2)$.
Then we have
$$
\frac{1}{d-1}\zeta(\rho^\theta)^{1-d}
\ge\int_{\rho^\theta}^{T/2} \frac{\zeta'(s)}{\zeta(s)^d}\,ds
\ge \zeta(\rho^\theta)^{-d}
\int_{\rho^\theta}^{T/2}
s^{\frac{N}{\theta}}f\left(C_*^{-1}s^{-\frac{N}{\theta}}\zeta(\rho^\theta)\right)\,ds.
$$
Recalling \eqref{eq:3.8} and setting $\eta:=C_*^{-1}\mu(B(z,\sigma)) s^{-\frac{N}{\theta}}$, 
we take large enough $C_*$ if necessary so that 
\begin{equation*}
\begin{split}
\mu(B(z,\sigma)) & \ge (d-1)\int_{\rho^\theta}^{T/2} s^{\frac{N}{\theta}}f\left(C_*^{-1}s^{-\frac{N}{\theta}}\zeta(\rho^\theta)\right)\,ds\\
 & =(d-1)\int_{\rho^\theta}^{T/2}
s^{\frac{N}{\theta}}f\left(C_*^{-1}\mu(B(z,\sigma))s^{-\frac{N}{\theta}}\right)\,ds\\
 & =\frac{(d-1)\theta}{N}C_*^{-p_\theta}\mu(B(z,\sigma))^{p_\theta}\int_{C_*^{-1}(T/2)^{-\frac{N}{\theta}}\mu(B(z,\sigma))}^{C_*^{-1}\rho^{-N}\mu(B(z,\sigma))}
\eta^{-p_\theta-1}f(\eta)\,d\eta\\
 & \ge \gamma^{-p_\theta-1}\mu(B(z,\sigma))^{p_\theta}\int_{\gamma^{-1}T^{-\frac{N}{\theta}}\mu(B(z,\sigma))}^{\gamma^{-1}\sigma^{-N}\mu(B(z,\sigma))}
\eta^{-p_\theta-1}f(\eta)\,d\eta
\end{split}
\end{equation*}
for all $z\in{\R}^N$ and $\sigma\in(0,T^{\frac{1}{\theta}})$, 
where $\gamma=2^{-\frac{N}{\theta}}C_*$. 
Here we used the relation $\sigma=2^{\frac{1}{\theta}}\rho\in(0,T^{\frac{1}{\theta}})$. 
Thus inequality~\eqref{eq:3.1} holds, and the proof of Theorem~\ref{Theorem:3.1} is complete.
\qed\vspace{5pt}
\begin{corollary}
\label{Corollary:3.1} 
Assume conditions~{\rm (F1)} and {\rm (F2)}.
Let $u$ satisfy
$$
\partial_t u+(-\Delta )^{\frac{\theta}{2}}u=F(u)
$$
in $Q_T$ for some $T>0$.
Then there exists a unique $\nu\in{\mathcal M}$ 
as the initial trace of $u$.
Furthermore,
\begin{itemize}
  \item[{\rm (i)}] 
  $u$ is a solution of problem~{\rm (P)} in $Q_T$ with $\mu=\nu$;
  \item[{\rm (ii)}]
  there exists $C=C(N,\theta,F)>0$ such that
  $$
  \sup_{z\in{{  \R}}^N}\mu(B(z,\sigma))\le
  \left\{
  \begin{array}{ll}
  C\sigma^{N-\frac{\theta}{p-1}}|\log\sigma|^{-\frac{q}{p-1}} &  \text{if}\quad p\neq p_\theta,\vspace{7pt}\\
  C|\log\sigma|^{-\frac{N(q+1)}{\theta}} &  \text{if}\quad p=p_\theta,\,\,\,q\neq-1,\vspace{7pt}\\
  C[\log|\log\sigma|]^{-\frac{N}{\theta}}  & \text{if}\quad p=p_\theta,\,\,\,q=-1,
  \end{array}
  \right.
  $$
  for all small enough $\sigma>0$. 
\end{itemize}
\end{corollary}
%
{\bf Proof.} 
The existence and  uniqueness of the initial trace of $u$ follows from Lemma~\ref{Lemma:2.4}. 
Let~$d\in (1,p)$, $R>0$, and $\kappa>0$. Set
$$
f(\tau):=
\left\{
\begin{array}{ll}
0 & \quad\mbox{for}\quad 0\le \tau<R,\vspace{5pt}\\
\displaystyle{\kappa\tau^d\int_R^\tau s^{-d}\left(\int_R^s \xi^{p-2}[\log(e+\xi)]^q\,d\xi\right)\,ds} & \quad\mbox{for}\quad \tau\ge R.
\end{array}
\right.
$$
By Lemma~\ref{Lemma:2.8}~(i) and (ii)
we see that $f$ is convex in $(0,\infty)$ and (f2) in Theorem~\ref{Theorem:3.1} holds. 
Furthermore, thanks to Lemma~\ref{Lemma:2.8}~(iii), 
taking small enough $\kappa>0$ and large enough $R>0$, 
by~(F2) we can ensure that $F(\tau)\ge f(\tau)$ in $[0,\infty)$ and consequently (f1) in Theorem~\ref{Theorem:3.1} also holds. 
In particular, we find $L\in(R,\infty)$ such that 
\begin{equation}
\label{eq:3.9}
F(\tau)\ge f(\tau)\ge C\tau^p(\log\tau)^q,\quad \tau\in(L,\infty). 
\end{equation}
By Theorem~\ref{Theorem:3.1} we also find $\gamma\ge 1$ such that
\begin{equation}
\label{eq:3.10}
\gamma^{p_\theta +1} m_\sigma(z)^{-\frac{\theta}{N}}\ge
\int_{\gamma^{-1}T^{-\frac{N}{\theta}}m_\sigma(z)}^{\gamma^{-1}\sigma^{-N}m_\sigma(z)}s^{-p_\theta-1}f(s)\,ds
\end{equation}
for all $z\in  \R^N$ and $\sigma\in(0,T^{\frac{1}{\theta}})$.

We show that 
\begin{equation}
\label{eq:3.11}
\sup_{z\in{{  \R}}^N}m_\sigma(z)<\infty\quad\mbox{for all $\sigma\in(0,T^{\frac{1}{\theta}})$}.
\end{equation}
For then 
by Remark~\ref{Remark:2.1}~(ii) we have
\begin{equation}
\label{eq:3.12}
\sup_{z\in{{  \R}}^N}\mu(B(z,1))\le C\sup_{z\in{{  \R}}^N}\mu(B(z,T^{\frac{1}{\theta}}/2))=C\sup_{z\in{{  \R}}^N}m_{T^{\frac{1}{\theta}}/2}(z)<\infty,
\end{equation}
and assertion~(i)  will follow from Lemma~\ref{Lemma:2.5}.

Suppose that $\sigma\in(0,T^{\frac{1}{\theta}})$ but \eqref{eq:3.11} does not hold. Then there exists
a sequence $\{z_n\}\subset\R^N$ such that $m_\sigma(z_n)\to\infty$ as $n\to\infty$. Consequently,
\begin{equation}
\label{eq:3.13}
\gamma^{-1}T^{-\frac{N}{\theta}}m_\sigma(z_n)\ge \max\{ L,2\}
\end{equation}
for all $n$ large enough.
By \eqref{eq:3.9}, \eqref{eq:3.10}, \eqref{eq:3.13}, and Lemma~\ref{Lemma:2.7} (i) (with $a=p-1$, $b=\theta/N$, and $c=q$),
 we obtain
\begin{equation*}
\begin{split}
& m_\sigma(z_n)^{-\frac{\theta}{N}}\\
& \ge C\gamma^{-p_\theta -1}\int_{\gamma^{-1}T^{-\frac{N}{\theta}}m_\sigma(z_n)}^{\gamma^{-1}\sigma^{-N}m_\sigma(z_n)}
s^{p-p_\theta-1}(\log s)^q\,ds\\
 & \ge C\left(\gamma^{-1}T^{-\frac{N}{\theta}}m_\sigma(z_n)\right)^{p-1}
 \left(\gamma^{-1}\sigma^{-N}m_\sigma(z_n)\right)^{-\frac{\theta}{N}}
\left(\log\left(\gamma^{-1}T^{-\frac{N}{\theta}}m_\sigma(z_n)\right)\right)^q\log\left(\frac{\sigma^{-N}}{T^{-\frac{N}{\theta}}}\right)\\
 &  =C\sigma^\theta T^{-\frac{N(p-1)}{\theta}}\log\left({\sigma^{-N}}{T^{\frac{N}{\theta}}}\right)m_\sigma(z_n)^{p-1-\frac{\theta}{N}}
 \left(\log\left(\gamma^{-1}T^{-\frac{N}{\theta}}m_\sigma(z_n)\right)\right)^q.
\end{split}
\end{equation*}
Hence
\begin{equation}
\label{eq:3.14}
1\ge
C\sigma^\theta T^{-\frac{N(p-1)}{\theta}}\log\left({\sigma^{-N}}{T^{\frac{N}{\theta}}}\right)m_\sigma(z_n)^{p-1}
 \left(\log\left(\gamma^{-1}T^{-\frac{N}{\theta}}m_\sigma(z_n)\right)\right)^q. 
\end{equation}
Letting $n\to\infty$ in \eqref{eq:3.14} yields a contradiction and  thus  \eqref{eq:3.11} holds. 

We now prove assertion~(ii).
Consider first the case where $p\not=p_\theta$. We show that there exist $C>0$ and  $\sigma_*>0$  such that 
\begin{equation}
\label{eq:3.15}
\sigma^{\frac{\theta}{p-1}-N}|\log\sigma|^{\frac{q}{p-1}}m_\sigma(z)\le C
\end{equation}
for all $z\in{  \R}^N$ and $\sigma\in(0,\sigma_*)$. Suppose, for contradiction, that there exist
 sequences $\{z_n\}\subset\R^N$ and $\{\sigma_n\}\subset(0,\infty)$ such that 
\begin{equation}
\label{eq:3.16}
\sigma_n\to 0\qquad\text{and}\qquad \sigma_n^{\frac{\theta}{p-1}-N}|\log\sigma_n|^{\frac{q}{p-1}}m_{\sigma_n}(z_n)\to\infty \qquad\text{as}\qquad  n\to\infty.
\end{equation}
Set $M_n:=m_{\sigma_n}(z_n)$. 
It follows from \eqref{eq:3.12} that 
\begin{equation}
\label{eq:3.17}
M_n\le \sup_{z\in{{  \R}}^N}\mu(B(z,1))\le C
\end{equation}
for all $n$ large enough. 
By \eqref{eq:3.16} we necessarily have
$$
 \sigma_n^{-N}M_n\to\infty \qquad\text{as}\qquad  n\to\infty,
$$
so that for $n$ large enough,
\begin{equation}
\label{eq:3.18}
 \gamma^{-1}(2\sigma_n)^{-N}M_n\ge \max\{ L,2\}\qquad\text{and}\qquad (2\sigma_n)^{-N}>T^{-\frac{N}{\theta}}. 
\end{equation}
Similar to the proof of part (i), it follows from \eqref{eq:3.9}, \eqref{eq:3.10}, and \eqref{eq:3.18} that
$$
\gamma^{p_\theta +1} \sigma_n^{-\theta}
\ge C\left(\sigma_n^{-N}M_n\right)^{\frac{\theta}{N}}\int_{\gamma^{-1}(2\sigma_n)^{-N}M_n}^{\gamma^{-1}\sigma_n^{-N}M_n}s^{p-p_\theta-1}(\log s)^q\,ds.
$$
Applying Lemma~\ref{Lemma:2.7} (i) (with $a=p-1$, $b=\theta/N$, and $c=q$), we obtain (for $n$ large enough) 
\begin{equation}
\label{eq:3.19}
 C\sigma_n^{-\theta}\ge \tau_n^{p-1}(\log(C\tau_n))^q,
\end{equation}
where $\tau_n:=\sigma_n^{-N}M_n$. 
For $n$ large enough, and rescaling  with $s_n=C\tau_n$ in \eqref{eq:3.19}, we can apply  Lemma~\ref{Lemma:2.6}
 (with $a=p-1>0$, $b=q$, and $c=0$) to obtain (after rescaling back to $\tau_n$) 
$$
\sigma_n^{-N}M_n=\tau_n
\le C\left(C\sigma_n^{-\theta}\right)^{\frac{1}{p-1}}\left(\log \left[C\sigma_n^{-\theta}\right]\right)^{-\frac{q}{p-1}}.
$$
Consequently, for such $n$,
$$
\sigma_n^{\frac{\theta}{p-1}-N}|\log\sigma_n|^{\frac{q}{p-1}}M_n\le C,
$$
contradicting \eqref{eq:3.16}. Thus \eqref{eq:3.15} holds, as required.

Now consider the case when $p=p_\theta$ and $q\not=-1$. 
We show that there exist $C>0$ and  $\sigma_*>0$  such that 
\begin{equation}
\label{eq:3.20}
|\log\sigma|^{\frac{N(q+1)}{\theta}}m_\sigma(z)\le C
\end{equation}
for all $z\in{  \R}^N$ and $\sigma\in(0,\sigma_*)$. Suppose, for contradiction, that there exist
sequences $\{z_n\}\subset\R^N$ and $\{\sigma_n\}\subset(0,\infty)$ such that 
\begin{equation}
\label{eq:3.21}
\sigma_n\to 0\qquad\text{and}\qquad |\log\sigma_n|^{\frac{N(q+1)}{\theta}}m_{\sigma_n}(z_n)
\to\infty \qquad\text{as}\qquad  n\to\infty.
\end{equation}
Set $M_n:=m_{\sigma_n}(z_n)$. Since  $\sigma^{-N/2}\ge |\log\sigma|^{\frac{N(q+1)}{\theta}}$ for all $\sigma >0$ small enough,
by \eqref{eq:3.21} we necessarily have
$$
 \sigma_n^{-\frac{N}{2}}M_n\to\infty \qquad\text{as}\qquad  n\to\infty,
$$
so that for $n$ large enough,
\begin{equation}
\label{eq:3.22}
 \gamma^{-1}\sigma_n^{-\frac{N}{2}}M_n\ge \max\{ L,2\}\qquad\text{and}\qquad \sigma_n^{-\frac{N}{2}}>T^{-\frac{N}{\theta}}. 
\end{equation}
Similar to the proof of part (i), it follows from \eqref{eq:3.9}, \eqref{eq:3.10}, and \eqref{eq:3.22} that
\begin{equation}
\label{eq:3.23}
\gamma^{p_\theta +1} 
  \ge CM_n^{\frac{\theta}{N}}\int_{\gamma^{-1}\sigma_n^{-\frac{N}{2}}M_n}^{\gamma^{-1}\sigma_n^{-N}M_n}s^{-1}(\log s)^q\,ds.
\end{equation}
Now set $c_q:=1/2$ if $q\ge 0$ and $c_q:=1$ if $q<0$. 
Then, by \eqref{eq:3.23} we have 
\begin{equation*}
\begin{split}
1 & \ge CM_n^{\frac{\theta}{N}}
 \left(\log(\gamma^{-1}\sigma_n^{-Nc_q} M_n)\right)^q\int_{\gamma^{-1}\sigma_n^{-\frac{N}{2}}M_n}
 ^{\gamma^{-1}\sigma_n^{-N}M_n}\tau^{-1}\,d\tau\\
 & = CM_n^{\frac{\theta}{N}}\left(\log \left(\gamma^{-1}\sigma_n^{-Nc_q}M_n\right)\right)^q\log\left(\sigma_n^{-\frac{N}{2}}\right),
\end{split}
\end{equation*}
so that
\begin{equation}
\label{eq:3.24}
\left(\sigma_n^{-Nc_q} M_n\right)^{\frac{\theta}{N}}
\left(\log \left(\gamma^{-1}\sigma_n^{-Nc_q}M_n\right)\right)^q
\le C\sigma_n^{-c_q\theta}|\log\sigma_n|^{-1}.
\end{equation}
Setting $\tau_n:=\gamma^{-1}\sigma_n^{-Nc_q}M_n$, \eqref{eq:3.24} can be written as 
\begin{equation}
\label{eq:3.25}
\tau_n^{\frac{\theta}{N}}\left(\log \tau_n\right)^q
\le C\sigma_n^{-c_q\theta}|\log\sigma_n|^{-1}.
\end{equation}
Applying Lemma~\ref{Lemma:2.6} to \eqref{eq:3.25} (with $a=\theta/N$, $b=q$, and $c=0$), 
then yields 
\begin{equation*}
\begin{split}
\sigma_n^{-Nc_q}M_n
& \le C\left(\sigma_n^{-c_q\theta}|\log\sigma_n|^{-1}\right)^{\frac{N}{\theta}}
\left(\log\left(C\sigma_n^{-c_q\theta}|\log\sigma_n|^{-1}\right)\right)^{-\frac{Nq}{\theta}}\\
&\le C\sigma_n^{-Nc_q}|\log\sigma_n|^{-\frac{N}{\theta}}
\left(C\log\left(\sigma_n^{-1}\right)\right)^{-\frac{Nq}{\theta}}\\
&\le C\sigma_n^{-Nc_q}|\log\sigma_n|^{-\frac{N(q+1)}{\theta}}
\end{split}
\end{equation*}
for all $n$ large enough. Consequently for such $n$,  
$$
|\log\sigma_n|^{\frac{N(q+1)}{\theta}}M_n\le C, 
$$
contradicting \eqref{eq:3.21}. Thus \eqref{eq:3.20} holds, as required.

Finally, consider the case when $p=p_\theta$ and $q=-1$. 
We show that there exist $C>0$ and  $\sigma_*>0$  such that 
\begin{equation}
\label{eq:3.26}
\left( \log|\log\sigma|\right)^{\frac{N}{\theta}}m_\sigma(z)\le C
\end{equation}
for all $z\in{  \R}^N$ and $\sigma\in(0,\sigma_*)$. Suppose, for contradiction, that there exist
sequences $\{z_n\}\subset\R^N$ and $\{\sigma_n\}\subset(0,\infty)$ such that 
\begin{equation}
\label{eq:3.27}
\sigma_n\to 0\qquad\text{and}\qquad \left( \log|\log\sigma_n|\right)^{\frac{N}{\theta}}m_{\sigma_n}(z_n)
\to\infty \qquad\text{as}\qquad  n\to\infty.
\end{equation}
Set $M_n:=m_{\sigma_n}(z_n)$. Since  $\sigma^{-N}\ge \left( \log|\log\sigma|\right)^{\frac{N}{\theta}}$ for all $\sigma >0$ small enough, 
by \eqref{eq:3.27} we necessarily have
$$
\sigma_n^{-N}M_n\to\infty \qquad\text{as}\qquad  n\to\infty.
$$
Then, combining \eqref{eq:3.17}, we find $L'>0$ such that 
\begin{equation}
\label{eq:3.28}
\max\{\gamma^{-1}T^{-\frac{N}{\theta}}M_n, L,2\}\le L'<\gamma^{-1}\sigma_n^{-N}M_n
\end{equation}
for all $n$ large enough. 
Once again, by \eqref{eq:3.9}, \eqref{eq:3.10}, and \eqref{eq:3.28} 
we have
\begin{equation*}
\begin{split}
 \sigma_n^{-\theta}
 & \ge C\gamma^{-p_\theta -1}\left(\sigma_n^{-N}M_n\right)^{\frac{\theta}{N}}
 \int_{L'}^{\gamma^{-1}\sigma_n^{-N}M_n}
\tau^{-1}(\log\tau)^{-1}\,d\tau\\
 & = C \tau_n^{\frac{\theta}{N}}\log\left(\frac{\log \tau_n}{\log L'}\right)\ge C\tau_n^{\frac{\theta}{N}}\log\log \tau_n
\end{split}
\end{equation*}
for all $n$ large enough, 
where $\tau_{n}:=\gamma^{-1}\sigma_n^{-N}M_n$.
By \eqref{eq:3.22} and  Lemma~\ref{Lemma:2.6} (with $a=\theta/N$, $b=0$, and $c=1$),
we have
$$
\gamma^{-1}\sigma_n^{-N}M_n=\tau_{n}\le \left( C\sigma_n^{-\theta}\right)^{\frac{N}{\theta}}
\left( \log\log \left[C\sigma_n^{-\theta}\right]\right)^{-\frac{N}{\theta}}
\le C\sigma_n^{-N}\left( \log|\log\sigma_n|\right)^{-\frac{N}{\theta}},
$$
so that 
$$
\left( \log|\log\sigma_n|\right)^{\frac{N}{\theta}}M_n\le C,
$$
contradicting \eqref{eq:3.27}. Hence \eqref{eq:3.26} holds, as required. 
The proof of Corollary~\ref{Corollary:3.1} is complete. 
\qed
\section{Sufficient Conditions  for Solvability }
\label{section:4}
In this section we establish sufficient conditions for the existence of  a supersolution, and consequently of a local-in-time solution of problem~(P), 
for three general classes of nonlinearity~$F$ (see Theorems~\ref{Theorem:4.1}, \ref{Theorem:4.2}, and \ref{Theorem:4.3}).
As corollaries, we  obtain the corresponding results when specializing to nonlinearities satisfying (F1) and (F2)
(Corollaries~\ref{Corollary:4.1}, \ref{Corollary:4.2}, and \ref{Corollary:4.3}). 
Indeed, for $F$ satisfying (F1) and (F2) the classification of initial data for which problem~(P) is locally solvable,
  separates naturally into the following three cases:
\begin{itemize}
  \item[(A):]
  $
  \mbox{either}\quad
  \mbox{{\rm (i)}}\quad\mbox{$1<p<p_\theta$ and $q\in{{\R}}$}\quad\mbox{or}\quad
  \mbox{{\rm (ii)}}\quad\mbox{$p=p_\theta$ and $q<-1$};
  $
   \item[(B):]
  $\mbox{$p=p_\theta$ and $q\ge-1$}$;
   \item[(C):]
  $p>p_\theta$.
\end{itemize}
\subsection{Sufficiency:  Case (A)}
\label{subsection:4.1}
We begin with nonlinearities $F$ which generalize case (A). 
\begin{theorem}
\label{Theorem:4.1}
Let $F$ be a nonnegative continuous function in $[0,\infty)$
and assume the following conditions:
\begin{itemize}
  \item[{\rm (A1)}]
  there exists $R\ge 0$ such that
  the function $(R,\infty)\ni \tau\mapsto \tau^{-1}F(\tau)$ is increasing;
  \item[{\rm (A2)}]
  $\displaystyle{\int_1^\infty \tau^{-p_\theta-1}F(\tau)\,d\tau<\infty}$.
\end{itemize}
If $\mu\in{\mathcal M}$ satisfies
\begin{equation}
\label{eq:4.1}
\sup_{x\in{{  \R}}^N}\mu(B(x,1))<\infty,
\end{equation}
then problem~{\rm (P)} possesses a solution $u$ in $Q_T$ for some $T>0$, with $u$ satisfying
$$
0\le u(x,t)\le 2[S(t)\mu](x)+R\le Ct^{-\frac{N}{\theta}}
$$
in $Q_T$ for some $C>0$.
\end{theorem}
{\bf Proof.}
Let $T\in(0,1)$ be chosen later. 
Set $w(x,t):=R+2[S(t)\mu](x)$.
It follows from Lemma~\ref{Lemma:2.1} and \eqref{eq:4.1} that
\begin{equation}
\label{eq:4.2}
R\le w(x,t)\le R+Ct^{-\frac{N}{\theta}}\sup_{z\in{{  \R}}^N}\mu(B(z,t^{\frac{1}{\theta}}))
\le R+Mt^{-\frac{N}{\theta}}\le 2Mt^{-\frac{N}{\theta}}
\end{equation}
for $0<t\le T$ and small enough $T$, where
$M:=C\displaystyle{\sup_{x\in{{  \R}}^N}}\mu(B(x,1))+1<\infty$.
Then, by (A1) and~\eqref{eq:4.2} we have
\begin{equation}
\label{eq:4.3}
\begin{split}
0\le\frac{F(w(x,t))}{w(x,t)}
\le (2M)^{-1}t^{\frac{N}{\theta}}F(2Mt^{-\frac{N}{\theta}}),
\quad (x,t)\in Q_T.
\end{split}
\end{equation}
Noting that
$$
S(t-s)w(s)=S(t-s)[R+2S(s)\mu]=R+2S(t)\mu=w(t),
$$
then by (A2) and \eqref{eq:4.3} we obtain
\begin{equation*}
\begin{split}
 & [S(t)\mu](x)+\int_0^t S(t-s)F(w(s))\,ds\\
 & \le \frac{1}{2}w(x,t)+\int_0^t \left\|\frac{F(w(s))}{w(s)}\right\|_{L^\infty({{  \R}}^N)} S(t-s)w(s)\,ds\\
 & \le \frac{1}{2}w(x,t)+(2M)^{-1}w(x,t)
 \int_0^t s^{\frac{N}{\theta}}F\left(2Ms^{-\frac{N}{\theta}}\right)\,ds\\
 & \le w(x,t)
  \left[\frac{1}{2}+CM^{\frac{\theta}{N}}
  \int_{2MT^{-\frac{N}{\theta}}}^\infty \tau^{-p_\theta-1}F(\tau)\,d\tau
  \right]\le w(x,t),
  \quad (x,t)\in Q_T,
\end{split}
\end{equation*}
for small enough $T$.  
This means that $w$ is a supersolution in $Q_T$
and the desired result follows from Lemma~\ref{Lemma:2.2} and \eqref{eq:4.2}. 
\qed\vspace{5pt}
\begin{corollary}
\label{Corollary:4.1}
Assume conditions~{\rm (F1)} and {\em (F2)} with
$$
\mbox{either}\quad
\mbox{{\rm (i)}}\quad\mbox{$1<p<p_\theta$ and $q\in{{  \R}}$}\qquad\mbox{or}\qquad
\mbox{{\rm (ii)}}\quad\mbox{$p=p_\theta$ and $q<-1$}.
$$
If $\mu\in{\mathcal M}$ satisfies
$$
\sup_{z\in{{\R}}^N}\mu(B(z,1))<\infty,
$$ 
then problem~{\rm (P)} possesses a solution~$u$ in $Q_T$ for some $T>0$,  
with $u$ satisfying 
$$
0\le u(x,t)\le 2[S(t)\mu](x)+R\le Ct^{-\frac{N}{\theta}}
$$
in $Q_T$ for some $R>0$ and $C>0$.
\end{corollary}
{\bf Proof.}
Set 
$$
g(\tau):=\tau\int_0^\tau s^{-1}\left(\int_0^s \xi^{p-2}[\log(e+\xi)]^q\,d\xi\right)\,ds, \qquad\tau\ge 0.
$$
It follows from Lemma~\ref{Lemma:2.8} (iii) (with $d=1$ and $R=0$) that $g(\tau)\asymp \tau^p[\log \tau]^q$ as $\tau\to\infty$. 
Hence, since either $1<p<p_\theta$, or $p=p_\theta$ and $q<-1$, 
we have  
\begin{equation}
\label{eq:4.4}
\int_1^\infty \tau^{-p_\theta-1}g(\tau)\,d\tau\le C \int_1^\infty \tau^{p-p_\theta-1}[\log \tau]^q\,d\tau<\infty.
\end{equation}
Let $\kappa>0$ and $L>0$. Set 
\begin{equation}
\label{eq:4.5}
f(\tau):=\kappa g(\tau)+L , \qquad\tau\ge 0.
\end{equation}
Clearly  $f(\tau)\asymp g(\tau)\asymp \tau^p[\log \tau]^q$ as $\tau\to\infty$ and so by (F1)--(F2)   we may 
choose  $\kappa$ and $L$ large enough such that 
\begin{equation}
\label{eq:4.6}
F(\tau)\le f(\tau),\quad \tau\ge 0.
\end{equation}
Now,
$$
\left(\frac{f(\tau)}{\tau}\right)'
=\kappa\tau^{-1}\left(\int_0^\tau \xi^{p-2}[\log(e+\xi)]^q\,d\xi\right)\,ds-L\tau^{-2}>0
$$
for all $\tau$ large enough ($\tau>R=R(\kappa , L)$). Hence $f$ satisfies hypothesis (A1) 
of Theorem~\ref{Theorem:4.1}. Furthermore, by \eqref{eq:4.4} and \eqref{eq:4.5}, $f$ also 
 satisfies hypothesis (A2) of Theorem~\ref{Theorem:4.1}.

Hence, by Theorem~\ref{Theorem:4.1}, there exists $T>0$ and a solution $v$  in $Q_T$  of problem~(P)  with $F$ replaced by $f$, with $v$ satisfying
$$
0\le v(x,t)\le 2[S(t)\mu](x)+R\le Ct^{-\frac{N}{\theta}}
$$
in $Q_T$ for some  $C>0$. 
This together with Lemma~\ref{Lemma:2.3} implies that problem~(P) possesses a solution $u$ in $Q_T$ such that
$$
0\le u(x,t)\le v(x,t)\le 2[S(t)\mu](x)+R\le Ct^{-\frac{N}{\theta}}
$$
in $Q_T$. Thus Corollary~\ref{Corollary:4.1} follows.
\qed
\subsection{Sufficiency: Case (B) }
We consider nonlinearities $F$ which generalize case (B).  
\begin{theorem}
\label{Theorem:4.2}
Let $\mu\in{\mathcal L}_0$ and let $F$ be an increasing, nonnegative continuous function in~$[0,\infty)$.
Assume that there exist $R>0$, $\alpha>0$, and positive functions $G\in C([R,\infty))$ and $H\in C^1([R,\infty))$
satisfying the following conditions {\rm (B1)-(B5)}:
\begin{itemize}
 \item[{\rm (B1)}] 
  $\tau^{-p_\theta}F(\tau)\asymp G(\tau)$ as $\tau\to\infty$;
 \item[{\rm (B2)}]
  {\rm (i)} for any $a\ge 1$ and $b>0$, 
  $G(a\tau^b)\asymp G(\tau)$ as $\tau\to\infty$. Furthermore, 
  {\rm (ii)} $\displaystyle{\lim_{\tau\to\infty}}\tau^{-\delta}G(\tau)=0$ for all $\delta>0$;
  \item[{\rm (B3)}]
  {\rm (i)} $H'(\tau)\asymp \tau^{-1}G(\tau)>0$ and 
  {\rm (ii)} $G(\tau H(\tau)^{-1})\asymp G(\tau)$ 
  as $\tau\to\infty$. Furthermore,  
  {\rm (iii)}~$\displaystyle{\lim_{\tau\to\infty}}H(\tau)=\infty$
   and {\rm (iv)} $\displaystyle{\lim_{\tau\to\infty}}\tau^{-\delta}H(\tau)=0$ 
  for all $\delta>0$;
  \item[{\rm (B4)}]
  there exists a strictly increasing and convex function $\Phi_\alpha$ in $[R,\infty)$ 
  such that 
  $$
  \Phi_\alpha^{-1}(\tau)=\tau H(\tau)^{-\alpha}
  $$ 
  for all $\tau\in[\Phi_\alpha(R),\infty)$;
  \item[{\rm (B5)}]
  there exists $\eta\in(0,\theta/N)$ such that the function 
  $P:(R,\infty)\ni\tau\mapsto \tau^\eta H(\tau)^{-\alpha} G(\tau)$ is increasing. 
\end{itemize}
Then there exists $\varepsilon>0$ such that if  $\mu$ satisfies
\begin{equation}
\label{eq:4.7}
\sup_{x\in{{  \R}}^N}
\Phi_\alpha^{-1}\biggr[\,\dashint_{B(x,\sigma)}\Phi_\alpha(\mu(y)+R)\,dy\biggr]\le\varepsilon\sigma^{-N}H(\sigma^{-1})^{-\frac{N}{\theta}}
\end{equation}
for all small enough $\sigma >0$,  
then problem~{\rm (P)} possesses a solution~$u$ in $Q_T$ for some {$T>0$}, 
with~$u$ satisfying
$$
0\le u(x,t)\le \Phi_\alpha^{-1}[S(t)\Phi_\alpha(\mu+C)]
\le Ct^{-\frac{N}{\theta}}H(t^{-1})^{-\frac{N}{\theta}}
$$
in $Q_T$ for some $C>0$.
\end{theorem}
We prepare a preliminary lemma.
\begin{lemma}
\label{Lemma:4.1}
Let $R>0$ and $\alpha>0$. Let $G$ and $H$ be positive functions in $[R,\infty)$
such that $G\in C([R,\infty))$ and $H\in C^1([R,\infty))$. 
Assume also that conditions {\rm (B2)-(i)}, {\rm (B3)-(i), (iii), (iv)}, and {\rm (B4)} in Theorem~{\rm\ref{Theorem:4.2}} hold. 
Then, for any $a>0$, $b>0$, and $c\in{\mathbb R}$, 
\begin{eqnarray}
\label{eq:4.8}
H(a\tau^b) & \asymp & H(\tau),\\
\label{eq:4.9}
 H(\tau^b H(\tau)^c) & \asymp & H(\tau),\\
\label{eq:4.10}
\Phi_\alpha(\tau) & \asymp & \tau H(\tau)^\alpha,
\end{eqnarray}
as $\tau\to\infty$. 
\end{lemma}
{\bf Proof.} 
We first prove \eqref{eq:4.8}. 
Consider the case where $a\ge 1$ and $b\ge 1$.  
By (B3)-(i) we see that $H$ is increasing for large enough $\tau$.  
Then we take large enough $R'\in(R,\infty)$ so that 
$a\tau^b\ge\tau\ge R'$ for $\tau\in[R',\infty)$, $(R'/a)^{1/b}\ge R$, and 
\begin{equation*}
\begin{split}
H(\tau)\le H(a\tau^b) & =\int_{R'}^{a\tau^b} H'(s)\,ds+H(R')
\le C\int_{R'}^{a\tau^b} s^{-1}G(s)\,ds+H(R')\\
 & =C\int^\tau_{(R'/a)^{1/b}} \xi^{-1}G(a\xi^b)\,d\xi+H(R')
\le C\int^\tau_R \xi^{-1}G(a\xi^b)\,d\xi+H(\tau)
\end{split}
\end{equation*}
for all $\tau\in[R',\infty)$, where $\xi=(s/a)^{1/b}$.
Then, by (B2)-(i) and (B3)-(i), (iii) we take large enough $R''\in(R',\infty)$ so that 
\begin{equation*}
\begin{split}
H(\tau) & \le H(a\tau^b)\le C\int^{R''}_R \xi^{-1}G(a\xi^b)\,d\xi+C\int^\tau_{R''} \xi^{-1}G(a\xi^b)\,d\xi +H(\tau)\\ 
 & \le C+C\int_{R''}^\tau \xi^{-1}G(\xi)\,d\xi+H(\tau)
 \le C+C\int_{R''}^\tau H'(\xi)\,d\xi+H(\tau)\le CH(\tau)+C\le CH(\tau)
\end{split}
\end{equation*}
for large enough $\tau$. Thus \eqref{eq:4.8} holds for $a\ge 1$ and $b\ge 1$. 
In particular, we have 
\begin{equation}
\label{eq:4.11}
H(\tau)\asymp H(a\tau)\asymp H(\tau^b)
\end{equation}
as $\tau\to\infty$ for $a\ge 1$ and $b\ge 1$. 
Then we see that 
\begin{equation}
\label{eq:4.12}
H(a^{-1}\tau)\asymp H(a\cdot a^{-1}\tau)=H(\tau),
\quad
H(\tau^{1/b})\asymp H((\tau^{1/b})^b)=H(\tau),
\end{equation}
as $\tau\to\infty$ for $a\ge 1$ and $b\ge 1$. 
By \eqref{eq:4.11} and \eqref{eq:4.12}, 
for any $a>0$ and $b>0$, we obtain 
$$
H(a\tau^b)\asymp H(\tau^b)\asymp H(\tau)
$$
as $\tau\to\infty$, and \eqref{eq:4.8} holds. 

Next we prove \eqref{eq:4.9}. Let $\delta>0$ be such that $b-\delta|c|>0$. 
By (B3)-(iii),~(iv) we see that 
$1\le H(\tau)\le \tau^{\delta}$ for large enough $\tau$. 
Since $H$ is increasing for large enough $\tau$, 
we have
$$
H(\tau^{b-|c|\delta})\le H(\tau^b H(\tau)^c)\le
H(\tau^{b+|c|\delta})
$$
as $\tau\to\infty$. 
This together with \eqref{eq:4.8} implies that 
$H(\tau^b H(\tau)^c)\asymp H(\tau)$  as $\tau\to\infty$, that is, 
\eqref{eq:4.9} holds.  
 
Furthermore, we observe from (B4) and \eqref{eq:4.9} (with $b=1$ and $c=\alpha$) that 
$$
\Phi_\alpha^{-1}(\tau H(\tau)^\alpha)=\tau H(\tau)^\alpha H(\tau H(\tau)^\alpha)^{-\alpha}\asymp\tau
$$
as $\tau\to\infty$.
Then we find $C\ge 1$ such that 
$$
C^{-1}\tau\le\Phi_\alpha^{-1}(\tau H(\tau)^\alpha)\le C\tau
$$
for large enough $\tau$, which together with (B4) implies that 
$$
\Phi_\alpha(C^{-1}\tau)\le \tau H(\tau)^\alpha\le\Phi_\alpha(C\tau)
$$
for large enough $\tau$. Then, by \eqref{eq:4.8} we see that
$$
\Phi_\alpha(\tau)\le C\tau H(C\tau)^\alpha\le C\tau H(\tau)^\alpha,
\qquad
\Phi_\alpha(\tau)\ge C^{-1}\tau H(C^{-1}\tau)^\alpha\ge C\tau H(\tau)^\alpha
$$
as $\tau\to\infty$, yielding \eqref{eq:4.10}. The proof of Lemma~\ref{Lemma:4.1} is complete.
\qed
\vspace{5pt}

{\bf Proof of Theorem~\ref{Theorem:4.2}.} 
Let $\varepsilon\in(0,1)$ and $L\in(R,\infty)$ be chosen later.
Set
$$
v(x,t):=[S(t)\Phi_\alpha(\mu+L)](x),\quad
w(x,t):=2\Phi_\alpha^{-1}(v(x,t)),\quad
\rho(\tau):=\tau^{-N}H(\tau^{-1})^{-\frac{N}{\theta}}.
$$
It follows from \eqref{eq:4.8} that
\begin{equation}
\label{eq:4.13}
\rho(t^{\frac{1}{\theta}})= t^{-\frac{N}{\theta}}H(t^{-\frac{1}{\theta}})^{-\frac{N}{\theta}}
\le Ct^{-\frac{N}{\theta}}H(t^{-1})^{-\frac{N}{\theta}}
\end{equation}
for all $t\in(0,T)$ and small enough $T$.
Furthermore, by (B3)-(iv) we see that $\rho(t^{\frac{1}{\theta}})\to \infty$ as~$t\to 0$. 
We apply Lemmas~\ref{Lemma:2.1} and \ref{Lemma:4.1} to obtain
\begin{align}
\Phi_\alpha(L)\le v(x,t) & \le  Ct^{-\frac{N}{\theta}}\sup_{z\in{{  \R}}^N}\int_{B(z,t^{\frac{1}{\theta}})}\Phi_\alpha(\mu(y)+L)\,dy
&   \text{[by Lemma~\ref{Lemma:2.1}]} \notag\\
 & \le  C\Phi_\alpha\left(\varepsilon\rho(t^{\frac{1}{\theta}})\right)
 \le C\varepsilon\rho(t^{\frac{1}{\theta}})H\left(\varepsilon\rho(t^{\frac{1}{\theta}})\right)^\alpha &  \text{[by \eqref{eq:4.7}, \eqref{eq:4.10}]} \notag\\
 & \le  C\varepsilon t^{-\frac{N}{\theta}}H(t^{-1})^{-\frac{N}{\theta}}
H\left(C\varepsilon t^{-\frac{N}{\theta}}H(t^{-1})^{-\frac{N}{\theta}}\right)^\alpha &   \text{[by (B3)-(i), \eqref{eq:4.13}]}
\label{eq:4.14}
\end{align}
in $Q_T$ for small enough $T$.
Since
\begin{equation}
\label{eq:4.15}
Ct^{-\frac{N}{\theta}} H(t^{-1})^{-\frac{N}{\theta}+\alpha}
\ge C\varepsilon t^{-\frac{N}{\theta}} H(t^{-1})^{-\frac{N}{\theta}+\alpha}
\ge t^{-\frac{N}{2\theta}} H(t^{-1})^{-\frac{N}{\theta}+\alpha}\to\infty 
\end{equation}
as $t\to 0$ (see (B3)-(iii),~(iv)), 
by (B3)-(i), \eqref{eq:4.8}, \eqref{eq:4.9}, and \eqref{eq:4.14} we have 
\begin{equation}
\label{eq:4.16}
\Phi_\alpha(L)\le v(x,t)
\le C\varepsilon t^{-\frac{N}{\theta}}H(t^{-1})^{-\frac{N}{\theta}}
H\left(Ct^{-\frac{N}{\theta}}H(t^{-1})^{-\frac{N}{\theta}}\right)^\alpha 
\le  C\varepsilon t^{-\frac{N}{\theta}} H(t^{-1})^{-\frac{N}{\theta}+\alpha}
\end{equation}
in $Q_T$ for small enough $T$.
By (B4) and \eqref{eq:4.16}
we have
\begin{equation}
\label{eq:4.17}
\begin{split}
2L\le w(x,t) & \le 2\Phi_\alpha^{-1}\left(C\varepsilon t^{-\frac{N}{\theta}} H(t^{-1})^{-\frac{N}{\theta}+\alpha}\right)\\
 & =C\varepsilon t^{-\frac{N}{\theta}} H(t^{-1})^{-\frac{N}{\theta}+\alpha}
 H\left(C\varepsilon t^{-\frac{N}{\theta}} H(t^{-1})^{-\frac{N}{\theta}+\alpha}\right)^{-\alpha}
\end{split}
\end{equation}
in $Q_T$  for small enough $T$. 
Since $H^{-\alpha}$ is monotone decreasing for large enough $\tau$, 
by \eqref{eq:4.15}
we have 
$$
H\left(C\varepsilon t^{-\frac{N}{\theta}} H(t^{-1})^{-\frac{N}{\theta}+\alpha}\right)^{-\alpha}
\le H\left(t^{-\frac{N}{2\theta}} H(t^{-1})^{-\frac{N}{\theta}+\alpha}\right)^{-\alpha}
$$
for all $t\in(0,T)$ and small enough $T$. 
This together with \eqref{eq:4.9} implies that 
\begin{equation}
\label{eq:4.18}
H\left(C\varepsilon t^{-\frac{N}{\theta}} H(t^{-1})^{-\frac{N}{\theta}+\alpha}\right)^{-\alpha}
\le CH(t^{-1})^{-\alpha}
\end{equation}
for all $t\in(0,T)$ and small enough $T$. 
By \eqref{eq:4.17} and \eqref{eq:4.18} we obtain
\begin{equation}
\label{eq:4.19}
2L\le w(x,t)  \le C\varepsilon t^{-\frac{N}{\theta}}H(t^{-1})^{-\frac{N}{\theta}}
\end{equation}
in $Q_T$  for small enough $T$. 
Then, taking large enough $L$ if necessary, 
by (B1) and \eqref{eq:4.10} we have
\begin{equation}
\label{eq:4.20}
\frac{F(w(x,t))}{v(x,s)}=\frac{F(w(x,t))}{\Phi_\alpha(w(x,t)/2)}
\le C\frac{w(x,t)^{p_\theta}G(w(x,t))}{w(x,t)H(w(x,t))^\alpha}
=Cw(x,t)^{\frac{\theta}{N}-\eta}P(w(x,t))
\end{equation}
in $Q_T$, where $P$ is as in (B5). 
Furthermore, 
by (B5) and \eqref{eq:4.19} we obtain 
\begin{equation}
\label{eq:4.21}
\begin{split}
 & P(w(x,t))\le P\left(  C\varepsilon t^{-\frac{N}{\theta}}H(t^{-1})^{-\frac{N}{\theta}}\right)
 \le P\left(  Ct^{-\frac{N}{\theta}}H(t^{-1})^{-\frac{N}{\theta}}\right)\\
 & =\left( Ct^{-\frac{N}{\theta}}H(t^{-1})^{-\frac{N}{\theta}}\right)^\eta
H\left( Ct^{-\frac{N}{\theta}}H(t^{-1})^{-\frac{N}{\theta}}\right)^{-\alpha}
G\left( Ct^{-\frac{N}{\theta}}H(t^{-1})^{-\frac{N}{\theta}}\right)
\end{split}
\end{equation}
in $Q_T$ for small enough $T$. 
On the other hand, 
by (B3)-(iv) we see that $t^{-1}H(t^{-1})^{-1}\to\infty$ as~$t\to 0$.
Then, by (B2)-(i) and (B3)-(ii) we see that
$$
G\left( Ct^{-\frac{N}{\theta}}H(t^{-1})^{-\frac{N}{\theta}}\right)
\le CG\left(t^{-1}H(t^{-1})^{-1}\right)\le CG(t^{-1})
$$
for all $t\in(0,T)$ and small enough $T$. 
This together with \eqref{eq:4.8}, \eqref{eq:4.9}, and \eqref{eq:4.21} implies that
\begin{equation}
\label{eq:4.22}
P(w(x,t))\le C\left( t^{-\frac{N}{\theta}}H(t^{-1})^{-\frac{N}{\theta}}\right)^\eta H(t^{-1})^{-\alpha}G(t^{-1})
\end{equation}
in $Q_T$. 
Since $0<\eta<\theta/N$ (see (B5)), 
by \eqref{eq:4.19}, \eqref{eq:4.20}, and \eqref{eq:4.22} we obtain 
\begin{equation}
\label{eq:4.23}
\begin{split}
\frac{F(w(x,t))}{v(x,s)}
 & \le C\left(C\varepsilon t^{-\frac{N}{\theta}}H(t^{-1})^{-\frac{N}{\theta}}\right)^{\frac{\theta}{N}-\eta}
\left( t^{-\frac{N}{\theta}}H(t^{-1})^{-\frac{N}{\theta}}\right)^\eta H(t^{-1})^{-\alpha}G(t^{-1})\\
 & \le C\varepsilon^{\frac{\theta}{N}-\eta} t^{-1}H(t^{-1})^{-1-\alpha}G(t^{-1})
\end{split}
\end{equation}
in $Q_T$.
Therefore, we deduce from (B3)-(i) and \eqref{eq:4.23} that
\begin{equation}
\label{eq:4.24}
\begin{split}
 \int_0^t \left\|\frac{F(w(s))}{v(s)}\right\|_\infty\,ds
 & \le C \varepsilon^{\frac{\theta}{N}-\eta}\int_0^t s^{-1}H(s^{-1})^{-1-\alpha}G(s^{-1})\,ds\\
 & \le C \varepsilon^{\frac{\theta}{N}-\eta}
\int_{t^{-1}}^\infty \tau^{-1}H(\tau)^{-1-\alpha}G(\tau)\,d\tau\\
 & \le C \varepsilon^{\frac{\theta}{N}-\eta}
\int_{t^{-1}}^\infty H(\tau)^{-1-\alpha}H'(\tau)\,d\tau\\
& \le C \varepsilon^{\frac{\theta}{N}-\eta}H(t^{-1})^{-\alpha}
\end{split}
\end{equation}
for all $t\in(0,T)$.
Similarly, by \eqref{eq:4.14} and Lemma~\ref{Lemma:4.1},
we have
\begin{equation}
\label{eq:4.25}
\begin{split}
 & \frac{[S(t)\Phi_\alpha(\mu+L)](x)}{w(x,t)}
 =\frac{v(x,t)}{2\Phi_\alpha^{-1}(v(x,t))}=\frac{1}{2}H(v(x,t))^\alpha
 \le CH(t^{-1})^\alpha
\end{split}
\end{equation}
in $Q_T$. Therefore,
taking small enough $\varepsilon\in(0,1)$, 
by \eqref{eq:4.24} and \eqref{eq:4.25}
we obtain
\begin{equation*}
\begin{split}
 & [S(t)\mu](x)+\int_0^t S(t-s)F(w(s))\,ds\\
 & \le \frac{1}{2}w(x,t)+\int_0^t \left\|\frac{F(w(s))}{v(s)}\right\|_\infty S(t-s)v(s)\,ds\\
 & =\frac{1}{2}w(x,t)+w(x,t)\frac{[S(t)\Phi_\alpha(\mu+L)](x)}{w(x,t)}\int_0^t \left\|\frac{F(w(s))}{v(s)}\right\|_\infty\,ds\\
 & \le w(x,t)\left[\frac{1}{2}+C\varepsilon^{\frac{\theta}{N}-\eta}\right]
 \le w(x,t)
\end{split}
\end{equation*}
in $Q_T$, 
where we have used the fact that 
$$
w(x,t)\ge 2\Phi^{-1}_\alpha\left( S(t)\Phi_\alpha (\mu)\right)\ge 2S(t)\mu
$$
by Jensen's inequality. 
Hence $w$ is a supersolution in $Q_T$ and 
Theorem~\ref{Theorem:4.2} now follows from Lemma~\ref{Lemma:2.2} and \eqref{eq:4.19}.
\qed
\begin{corollary}
\label{Corollary:4.2}
Let $\mu\in{\mathcal L}_0$ and assume conditions~{\rm (F1)} and {\em (F2)} hold with $p=p_\theta$ and $q\ge -1$.
Let $\alpha>0$ and set
$$
\begin{array}{ll}
h(\tau) & :=
\left\{\begin{array}{ll}
(\log(e+\tau))^{q+1} & \mbox{if}\quad q>-1,\vspace{5pt}\\
\log(e+\log(e+\tau)) & \mbox{if}\quad q=-1,
\end{array}
\right.\vspace{5pt}\\
\psi_\alpha^\pm(\tau) & :=\tau h(\tau)^{\pm\alpha},
\end{array}
$$
for $\tau\in(0,\infty)$. 
Then there exists $\varepsilon>0$ such that 
if $\mu$ satisfies
\begin{equation}
\label{eq:4.26}
\sup_{x\in{{  \R}}^N}
\psi_\alpha^-\biggr[\,\dashint_{B(x,\sigma)}\psi_\alpha^+(\mu)\,dy\biggr]\le\varepsilon\sigma^{-N}h(\sigma^{-1})^{-\frac{N}{\theta}}
\end{equation}
for small enough $\sigma>0$, 
then problem~{\rm (P)} possesses a solution~$u$ in $Q_T$ for some $T>0$,
with $u$ satisfying 
$$
0\le u(x,t)\le
\left\{
\begin{array}{ll}
Ct^{-\frac{N}{\theta}}|\log t|^{-\frac{N(q+1)}{\theta}} & \mbox{if}\quad q>-1,\vspace{5pt}\\
Ct^{-\frac{N}{\theta}}[\log|\log t|]^{-\frac{N}{\theta}} & \mbox{if}\quad q=-1
\end{array}
\right.
$$
in $Q_T$ for some $C>0$.
\end{corollary}
{\bf Proof.}
Let $\alpha>1$. Set 
$$
G(\tau):=(\log \tau)^q,\quad
H(\tau):=
\left\{
\begin{array}{ll}
(\log \tau)^{q+1} & \mbox{if}\quad q>-1,\vspace{5pt}\\
\log(\log \tau) & \mbox{if}\quad q=-1,
\end{array}
\right.
\quad
\Psi_\alpha(\tau):=\tau H(\tau)^{-\alpha}.
$$
Then, for any $a\ge 1$ and $b>0$, we have 
\begin{align*}
G(a\tau^b) & =(\log a+b\log \tau)^q=b^q(\log\tau)^q(1+o(1))\\
 & \asymp (\log\tau)^q=G(\tau)
 \quad\mbox{for any $a\ge 1$ and $b>0$},\\
G(\tau) & =o(\tau^\delta)\quad\mbox{for any $\delta>0$}, 
\end{align*}
as $\tau\to\infty$. Thus condition~(B2) holds.  
Let $\Phi_\alpha$ be the inverse function of $\Psi_\alpha$, 
and we show that conditions~(B3)--(B5) in Theorem~\ref{Theorem:4.2} hold. 

Consider the case of $q>-1$. 
Then 
\begin{equation*}
\begin{split}
H'(\tau) & =(q+1)\tau^{-1}(\log\tau)^q,\\
\Psi_\alpha'(\tau) & =(\log\tau)^{-\alpha(q+1)}-\alpha(q+1)(\log\tau)^{-\alpha(q+1)-1}=(\log\tau)^{-\alpha(q+1)}(1+o(1))>0,\\
\Psi_\alpha''(\tau) & =-\alpha(q+1)\tau^{-1}(\log\tau)^{-\alpha(q+1)-1}+\alpha(q+1)(\alpha(q+1)+1)\tau^{-1}(\log\tau)^{-\alpha(q+1)-2}\\
 & =-\alpha(q+1)\tau^{-1}(\log\tau)^{-\alpha(q+1)-1}(1+o(1))<0,
\end{split}
\end{equation*}
as $\tau\to\infty$. 
We see that 
\begin{align*}
 & \tau^{-1}G(\tau)=\tau^{-1}(\log\tau)^q\asymp H'(\tau)\quad\mbox{as $\tau\to\infty$},\\
 & G(\tau H(\tau)^{-1})=\left[\log(\tau(\log\tau)^{-(q+1)})\right]^q\asymp (\log\tau)^q=G(\tau)\quad\mbox{as $\tau\to\infty$},\\
 & \lim_{\tau\to\infty}H(\tau)=\infty,\qquad H(\tau)=o(\tau^\delta)\quad \mbox{as $\tau\to\infty$ for any $\delta>0$}.
\end{align*}
Thus condition~(B3) holds. Furthermore, we observe that $\Psi_\alpha$ is strictly increasing and concave for large enough $\tau$, 
that is, the inverse function $\Phi_\alpha$ of $\Psi_\alpha^{-1}$ exists and it is strictly increasing and convex for large enough $\tau$. 
Then conditions~(B4) holds. 
In addition, 
for any $\eta>0$, setting 
$$
P(\tau)=\tau^\eta H(\tau)^{-\alpha}G(\tau)=\tau^\eta(\log\tau)^{-\alpha(q+1)+q},
$$
by Lemma~\ref{Lemma:2.6} we see that $P'(\tau)>0$ for large enough $\tau$. 
This implies that condition~(B5) also holds. 
Thus conditions~(B3)--(B5) hold in the case of $q>-1$. 
 
Consider the case of $q=-1$. 
It follows that
\begin{equation*}
\begin{split}
H'(\tau) & =\tau^{-1}(\log\tau)^{-1},\\
\Psi_\alpha'(\tau) & =(\log(\log\tau))^{-\alpha}-\alpha(\log\tau)^{-1}(\log(\log\tau))^{-\alpha-1}
=(\log(\log\tau))^{-\alpha}(1+o(1))>0,\\
\Psi_\alpha''(\tau) & =-\alpha\tau^{-1}(\log\tau)^{-1}(\log(\log\tau))^{-\alpha-1}
+\alpha\tau^{-1}(\log\tau)^{-2}(\log(\log\tau))^{-\alpha-1}\\
 & \qquad\quad
+\alpha(\alpha+1)\tau^{-1}(\log\tau)^{-2}(\log(\log\tau))^{-\alpha-2}\\
 & =-\alpha\tau^{-1}(\log\tau)^{-1}(\log(\log\tau))^{-\alpha-1}(1+o(1))<0,
\end{split}
\end{equation*}
as $\tau\to\infty$. 
Similarly to the case of $q>-1$, 
we have 
\begin{align*}
 & \tau^{-1}G(\tau)=\tau^{-1}(\log\tau)^{-1}=H'(\tau)\quad\mbox{as $\tau\to\infty$},\\
 & G(\tau H(\tau)^{-1})=\left[\log(\tau(\log(\log\tau))^{-1})\right]^{-1}\asymp (\log\tau)^{-1}=G(\tau)\quad\mbox{as $\tau\to\infty$},\\
 & \lim_{\tau\to\infty}H(\tau)=\infty,\qquad H(\tau)=o(\tau^\delta)\quad \mbox{as $\tau\to\infty$ for any $\delta>0$}.
\end{align*}
Thus condition~(B3) holds. 
Furthermore, we see that $\Psi_\alpha$ is strictly increasing and concave for large enough $\tau$, 
that is, the inverse function $\Psi_\alpha^{-1}$ exists and it is strictly increasing and convex for large enough $\tau$. 
Then conditions~(B4) holds. 
In addition, 
for any $\eta>0$, 
setting 
$$
P(\tau)=\tau^\eta H(\tau)^{-\alpha}G(\tau)=\tau^\eta(\log(\log\tau))^{-\alpha}(\log\tau)^{-1},
$$
by Lemma~\ref{Lemma:2.6} we see that $P'(\tau)>0$ for large enough $\tau$. 
This implies that condition~(B5) also holds. 
Thus conditions~(B3)--(B5) hold in the case of $q=-1$. 

Assume \eqref{eq:4.26}. 
By Lemma~\ref{Lemma:2.6} (with $a=1$, $b=-\alpha(q+1)$, and $c=0$ for $q>-1$ 
and with $a=1$, $b=0$, and $c=-\alpha$ for $q=-1$) 
we have 
$$
\Phi_\alpha(\tau)=\Psi_\alpha^{-1}(\tau)\asymp 
\left\{
\begin{array}{ll}
\tau(\log\tau)^{\alpha(q+1)} & \mbox{for}\quad q>-1,\vspace{5pt}\\
\tau(\log(\log\tau))^\alpha  & \mbox{for}\quad q=-1,
\end{array}
\right.
$$
as $\tau\to\infty$. 
Since 
$$
\Phi_\alpha^{-1}(\tau)=\Psi_\alpha(\tau)\le C\psi_\alpha^-(\tau),
\qquad
\Phi_\alpha(\tau)=\Psi_\alpha^{-1}(\tau)\le C\psi_\alpha^{+}(\tau),
$$
for large enough $\tau$,  
taking large enough $R>0$ if necessary, we see that  
\begin{equation}
\label{eq:4.27}
\Phi_\alpha^{-1}\left[\,\dashint_{B(x,\sigma)}\Phi_\alpha(\mu(y)+R)\,dy\right]
\le C\psi_\alpha^-\left[\,\dashint_{B(x,\sigma)}C\psi_\alpha^+(\mu(y)+R)\,dy\right].
\end{equation}
Furthermore, we see that 
$$
\psi_\alpha^+(\tau+R)\le C\psi_\alpha^+(\tau)+C,
\quad
\psi_\alpha^-(C\tau+C)\le C\psi_\alpha^-(\tau)+C
$$
for $\tau>0$. 
Then, by \eqref{eq:4.26} and \eqref{eq:4.27} we see that
\begin{equation*}
\begin{split}
 & \Phi_\alpha^{-1}\left[\,\dashint_{B(x,\sigma)}\Phi_\alpha(\mu(y)+R)\,dy\right]
\le C\psi_\alpha^-\left[\,\dashint_{B(x,\sigma)}C\psi_\alpha^+(\mu(y))\,dy+C\right]\\
 & \le C\psi_\alpha^-\left[\,\dashint_{B(x,\sigma)}\psi_\alpha^+(\mu(y))\,dy\right]+C\\
 & \le C\varepsilon\sigma^{-N}h(\sigma^{-1})^{-\frac{N}{\theta}}+C
 \le C\varepsilon\sigma^{-N}H(\sigma^{-1})^{-\frac{N}{\theta}}
\end{split}
\end{equation*}
for all small enough $\sigma>0$. 

Let $f$ be as in \eqref{eq:4.5}. 
Since $f(\tau)\asymp\tau^{p_\theta}(\log\tau)^q$ as $\tau\to\infty$, 
condition~(B1) holds with $F$ replaced by $f$. 
We deduce from Theorem~\ref{Theorem:4.2} that  
problem~(P) with $F$ replaced by $f$
possesses a solution $v$ in $Q_T$ for some $T>0$ such that
$$
0\le v(x,t)\le
\left\{
\begin{array}{ll}
Ct^{-\frac{N}{\theta}}|\log t|^{-\frac{N(q+1)}{\theta}} & \mbox{if}\quad q>-1,\vspace{5pt}\\
Ct^{-\frac{N}{\theta}}[\log|\log t|]^{-\frac{N}{\theta}} & \mbox{if}\quad q=-1,
\end{array}
\right.
$$
for all $(x,t)\in Q_T$.
This together with $f(\tau)\ge F(\tau )$ (by \eqref{eq:4.6}) and Lemma~\ref{Lemma:2.3} 
implies that problem~(P) possesses a solution $u$ in $Q_T$ such that
$$
0\le u(x,t)\le v(x,t)\le
\left\{
\begin{array}{ll}
Ct^{-\frac{N}{\theta}}|\log t|^{-\frac{N(q+1)}{\theta}} & \mbox{if}\quad q>-1,\vspace{5pt}\\
Ct^{-\frac{N}{\theta}}[\log|\log t|]^{-\frac{N}{\theta}} & \mbox{if}\quad q=-1,
\end{array}
\right.
$$
for all $(x,t)\in Q_T$.
Thus Corollary~\ref{Corollary:4.2} follows.
\qed

\subsection{Sufficiency:   Case (C)}
\label{subsection:4.3}
 In this section we consider nonlinearities $F$ which generalize case~(C). 
\begin{theorem}
\label{Theorem:4.3} 
Let $\mu\in{\mathcal L}_0$ and 
let $F$ be an increasing, nonnegative continuous function in~$[0,\infty)$ such that
\begin{itemize}
\item[{\rm (C1)}]
there exist $R\ge 0$ and $d>1$ such that
the function $(R,\infty)\ni \tau\mapsto \tau^{-d}F(\tau)\in(0,\infty)$ is increasing.
\end{itemize}
Furthermore,
assume that there exists a continuous function $G$ in $[R,\infty)$
satisfying the following conditions:
\begin{itemize}
\item[{\rm (C2)}]
there exists $p\in[d,d+1)$ such that $G(\tau)\succeq \tau^{-p}F(\tau)>0$ as $\tau \to \infty$;
\item[{\rm (C3)}]
for any $a\ge 1$, $b>0$, and $c\in{{  \R}}$,
$G(a\tau^b G(\tau)^c)\asymp G(\tau)$ as $\tau\to\infty$;
\item[{\rm (C4)}]
there exists $\delta\in(0,1)$ such that the function
$(R,\infty)\ni\tau\mapsto\tau^{-\delta}G(\tau)$ is decreasing.
\end{itemize}
Let $\alpha>1$.
Then there exists $\varepsilon>0$ such that
if $\mu$ satisfies
\begin{equation}
\label{eq:4.28}
\sup_{x\in{{  \R}}^N}\left[\,\dashint_{B(x,\sigma)}
\mu(y)^\alpha\,dy\,\right]^{\frac{1}{\alpha}}\le\varepsilon\sigma^{-\frac{\theta}{p-1}}G(\sigma^{-1})^{-\frac{1}{p-1}}
\end{equation}
for small enough $\sigma >0$, 
then problem~{\rm (P)} possesses a solution~$u$ in $Q_T$ for some {$T>0$}, 
with $u$ satisfying
$$
0\le u(x,t)\le 2[S(t)\mu^\alpha](x)^{\frac{1}{\alpha}}+R\le Ct^{-\frac{1}{p-1}} G(t^{-1})^{-\frac{1}{p-1}}
$$
in $Q_T$ for some $C>0$.
\end{theorem}
{\bf Proof.}
Let $\varepsilon\in(0,1)$ be chosen later, and assume \eqref{eq:4.28}. 
Without loss of generality we may assume  that $\alpha\in(1,d)$. 
 Indeed, if $\alpha\ge d$ and \eqref{eq:4.28} holds, then for any $\alpha'\in(1,d)$ we can write $\mu^{\alpha}=(\mu^{\alpha'})^\frac{\alpha}{\alpha'}$ 
 and apply Jensen's inequality to give
$$
\sup_{x\in{{  \R}}^N}\left[\,\dashint_{B(x,\sigma)}
\mu(y)^{\alpha'}\,dy\,\right]^{\frac{1}{\alpha'}}
\le\sup_{x\in{{  \R}}^N}
\left[\,\dashint_{B(x,\sigma)}
\mu(y)^\alpha\,dy\,\right]^{\frac{1}{\alpha}}
\le \varepsilon\sigma^{-\frac{\theta}{p-1}}G(\sigma^{-1})^{-\frac{1}{p-1}}
$$
for small enough $\sigma>0$. Consequently \eqref{eq:4.28} also holds for $\alpha'\in(1,d)$.

Set
\begin{equation}
\label{eq:4.29}
w(x,t):=2[S(t)\mu^\alpha](x)^{\frac{1}{\alpha}}+R.
\end{equation}
It follows from (C3), Lemma~\ref{Lemma:2.1}, and \eqref{eq:4.28} that
\begin{equation*}
\begin{split}
0 & \le [S(t)\mu^\alpha](x) \le Ct^{-\frac{N}{\theta}}\sup_{z\in{{  \R}}^N}\int_{B(z,t^{\frac{1}{\theta}})}\mu(y)^\alpha\,dy\\
 & \le C\varepsilon^\alpha t^{-\frac{\alpha}{p-1}}G(t^{-\frac{1}{\theta}})^{-\frac{\alpha}{p-1}}
 \le C\varepsilon^\alpha t^{-\frac{\alpha}{p-1}}G(t^{-1})^{-\frac{\alpha}{p-1}}
\end{split}
\end{equation*}
in $Q_T$ for small enough $T$. 
On the other hand, by (C1) and (C2) we see that
$$
\lim_{\tau\to\infty}\tau G(\tau)\ge C\lim_{\tau\to\infty}\tau^{1-p}F(\tau)=C\lim_{\tau\to\infty}\tau^{d+1-p}\tau^{-d}F(\tau)=\infty,
$$
since $p<d+1$. These imply that
\begin{align}
\label{eq:4.30}
R\le w(x,t) \le R+C \varepsilon t^{-\frac{1}{p-1}} G(t^{-1})^{-\frac{1}{p-1}}
 & \le C \varepsilon t^{-\frac{1}{p-1}} G(t^{-1})^{-\frac{1}{p-1}}\\
\label{eq:4.31}
  & \le Ct^{-\frac{1}{p-1}} G(t^{-1})^{-\frac{1}{p-1}}
\end{align}
in $Q_T$. 
Since $1<\alpha<d$,
by (C1)--(C3), \eqref{eq:4.30}, and \eqref{eq:4.31}
we obtain
\begin{equation}
\label{eq:4.32}
\begin{split}
&\frac{F(w(x,t))}{w(x,t)^\alpha}
=w(x,t)^{{d-\alpha}}\frac{F(w(x,t))}{w(x,t)^d}\\
 & \le C\left[C\varepsilon t^{-\frac{1}{p-1}}G(t^{-1})^{-\frac{1}{p-1}}\right]^{{d-\alpha}}
 \left[Ct^{-\frac{1}{p-1}}G(t^{-1})^{-\frac{1}{p-1}}\right]^{-d}
 F\left(Ct^{-\frac{1}{p-1}}G(t^{-1})^{-\frac{1}{p-1}}\right)\\
 & \le C\left[C\varepsilon t^{-\frac{1}{p-1}}G(t^{-1})^{-\frac{1}{p-1}}\right]^{{d-\alpha}}
 \left[Ct^{-\frac{1}{p-1}}G(t^{-1})^{-\frac{1}{p-1}}\right]^{p-d}
 G\left(Ct^{-\frac{1}{p-1}}G(t^{-1})^{-\frac{1}{p-1}}\right)\\
 & \le C\varepsilon^{d-\alpha}t^{-\frac{p-\alpha}{p-1}}G(t^{-1})^{\frac{\alpha-1}{p-1}}
\end{split}
\end{equation}
in $Q_T$ for small enough $T$. 
Similarly, by \eqref{eq:4.31} we have
\begin{equation}
\label{eq:4.33}
\begin{split}
{w(x,t)^{\alpha-1}}
\le Ct^{-\frac{\alpha-1}{p-1}}G(t^{-1})^{-\frac{\alpha-1}{p-1}}
\end{split}
\end{equation}
in $Q_T$. 
On the other hand, 
by (C4) we see that
\begin{equation}
\label{eq:4.34}
\begin{split}
 & \int_0^t s^{-\frac{p-\alpha}{p-1}}G(s^{-1})^{\frac{\alpha-1}{p-1}}\,ds
 =\int_{t^{-1}}^\infty \tau^{\frac{p-\alpha}{p-1}-2}G(\tau)^{\frac{\alpha-1}{p-1}}\,d\tau\\
 & =\int_{t^{-1}}^\infty \tau^{-1-(1-\delta)\frac{\alpha-1}{p-1}}
[\tau^{-\delta}G(\tau)]^{\frac{\alpha-1}{p-1}}\,d\tau\\
 & \le C[t^\delta G(t^{-1})]^{\frac{\alpha-1}{p-1}}t^{(1-\delta)\frac{\alpha-1}{p-1}}
=Ct^{\frac{\alpha-1}{p-1}}G(t^{-1})^{\frac{\alpha-1}{p-1}}
\end{split}
\end{equation}
for all $t\in (0,T)$ and small enough $T$. 
Therefore,
taking small enough $\varepsilon$,
by Jensen's inequalities, \eqref{eq:4.32}, \eqref{eq:4.33}, and \eqref{eq:4.34}
we obtain
\begin{equation*}
\begin{split}
 & [S(t)\mu](x)+\int_0^t S(t-s)F(w(s))\,ds\\
 & \le [S(t)\mu^\alpha](x)^{\frac{1}{\alpha}}+C\int_0^t \left\|\frac{F(w(s))}{w(s)^\alpha}\right\|_{L^\infty({{  \R}}^N)}
 S(t-s)[S(s)\mu^\alpha+R^\alpha]\,ds\\
 & \le \frac{1}{2}w(x,t)+C\varepsilon^{d-\alpha}[S(t)\mu^\alpha+R^\alpha]
 \int_0^t  s^{-\frac{p-\alpha}{p-1}}G(s^{-1})^{\frac{\alpha-1}{p-1}}\,ds\\
  & \le \frac{1}{2}w(x,t)+C\varepsilon^{d-\alpha}w(x,t)^\alpha
 \int_0^t  s^{-\frac{p-\alpha}{p-1}}G(s^{-1})^{\frac{\alpha-1}{p-1}}\,ds\\
 & \le \frac{1}{2}w(x,t)
 +C\varepsilon^{d-\alpha}\left\|w(t)^{\alpha-1}\right\|_{L^\infty({{  \R}}^N)} t^{\frac{\alpha-1}{p-1}}G(t^{-1})^{\frac{\alpha-1}{p-1}} w(x,t)\\
 & \le w(x,t)\left[\frac{1}{2}+C\varepsilon^{d-\alpha}\right]\le w(x,t)
\end{split}
\end{equation*}
in $Q_T$. 
Hence $w$ is a supersolution in $Q_T$ and Theorem~\ref{Theorem:4.3} now follows from Lemma~\ref{Lemma:2.2}, \eqref{eq:4.29}, and \eqref{eq:4.30}.
\qed\vspace{3pt}
\begin{corollary}
\label{Corollary:4.3}
Let $\mu\in{\mathcal L}_0$ and assume conditions~{\rm (F1)} and {\em (F2)} hold.
For any $\alpha>1$,
there exists $\varepsilon>0$ such that if $\mu$ satisfies
\begin{equation}
\label{eq:4.35}
\sup_{x\in{{  \R}}^N}\left[\,\dashint_{B(x,\sigma)}
\mu(y)^\alpha\,dy\,\right]^{\frac{1}{\alpha}}\le\varepsilon\sigma^{-\frac{\theta}{p-1}}
|\log\sigma|^{-\frac{q}{p-1}}
\end{equation}
for all small enough $\sigma>0$, 
then problem~{\rm (P)}  possesses a solution~$u$ in $Q_T$ for some $T>0$,
with~$u$ satisfying
$$
0\le u(x,t)\le 2[S(t)\mu^\alpha](x)^{\frac{1}{\alpha}}+R\le Ct^{-\frac{1}{p-1}}|\log t|^{-\frac{q}{p-1}}
$$
in $Q_T$ for some $R, C>0$.
\end{corollary}
{\bf Proof.}
Let $d\in(1,p)$ with $d>p-1$. 
Let $\kappa$, $L>0$, and set
\begin{equation}
\label{eq:4.36}
f(\tau):=\kappa\tau^d\int_0^\tau s^{-d}\left(\int_0^s \xi^{p-2}[\log(e+\xi)]^q\,d\xi\right)\,ds+L
\end{equation}
for $\tau\in  (0,\infty)$.
It follows from Lemma~\ref{Lemma:2.8}~(iii) that 
\begin{equation}
\label{eq:4.37}
\tau^d\int_0^\tau s^{-d}\left(\int_0^s \xi^{p-2}[\log(e+\xi)]^q\,d\xi\right)\,ds\asymp \tau^p[\log \tau]^q
\end{equation}
as $\tau\to\infty$. 
We take large enough $\kappa$ and $L$ so that $F(\tau)\le f(\tau)$ in $[0,\infty)$. 
On the other hand, since
$$
\left(\frac{f(\tau)}{\tau^d}\right)'=\tau^{-d}
\left[
\kappa\left(\int_0^\tau \xi^{p-2}[\log(e+\xi)]^q\,d\xi\right)\,ds
-Ld\tau^{-1}\right]>0
$$
for large enough $\tau$, 
condition~(C1) in Theorem~\ref{Theorem:4.3} holds in $(R,\infty)$ with $F$ replaced by $f$ for some $R>0$. 

Taking large enough $R$ if necessary and 
setting $G(\tau):=(\log\tau)^q$ for $\tau\in(R,\infty)$, 
we see that the function $(R,\infty)\ni\tau\mapsto \tau^{-\frac{1}{2}}G(\tau)$ is decreasing (i.e. $\delta =1/2$ in (C4)). 
By \eqref{eq:4.36} and \eqref{eq:4.37} we find $C>0$ such that 
$$
\tau^{-p}f(\tau)\le CG(\tau)
$$
for all $\tau\in(R,\infty)$. 
Then conditions~(C2)--(C4) in Theorem~\ref{Theorem:4.3} hold with $F$ replaced by $f$. 
Therefore, by Theorem~\ref{Theorem:4.3} 
there exists $\varepsilon>0$ such that
if $\mu$ satisfies \eqref{eq:4.35}, then problem~(P) with $F$ replaced by $f$
possesses a solution $v$ in $Q_T$ for some $T>0$ such that
$$
0\le v(x,t)\le 2[S(t)\mu^\alpha](x)^{\frac{1}{\alpha}}+R\le Ct^{-\frac{1}{p-1}}|\log t|^{-\frac{q}{p-1}}
$$
in $Q_T$, for some $C>0$. 
This together with Lemma~\ref{Lemma:2.3} 
implies that problem~(P) possesses a solution $u$ in $Q_T$ such that
$$
0\le u(x,t)\le v(x,t)\le 2[S(t)\mu^\alpha](x)^{\frac{1}{\alpha}}+R\le Ct^{-\frac{1}{p-1}}|\log t|^{-\frac{q}{p-1}}
$$
in $Q_T$. Thus Corollary~\ref{Corollary:4.3} follows.
\qed
%
\subsection{A special case:  Dirac measure as initial data}
\label{subsection:4.4}
Here we provide a necessary and sufficient condition on the nonlinearity $F$
for the solvability of problem~(P) in the special case when $\mu=\delta_y$, the Dirac measure in ${{  \R}}^N$ based at point $y$. 
This problem was considered in \cite{BF}
for the opposite sign pure power law case $F(u)=-u^p$, i.e. dissipative $F$. 
\begin{corollary}
\label{Corollary:4.4}
Suppose $F$ satisfies
\begin{itemize}
  \item[{\rm (D1)}]
  $F$ is nonnegative and locally Lipschitz continuous in $[0,\infty)$;
  \item[{\rm (D2)}]
  there exist $R>0$ and $d>1$ such that
  \begin{itemize}
  \item[{\rm (i)}]
  the function $(R,\infty)\ni \tau\mapsto \tau^{-d}F(\tau)\in(0,\infty)$ is increasing;
  \item[{\rm (ii)}]
  $F$ is convex in $(R,\infty)$.
  \end{itemize}
\end{itemize}
Let $y\in\R^N$. Then problem~{\rm (P)} possesses a local-in-time solution with $\mu=\delta_y$ if and only if
\begin{equation}
\label{eq:4.38}
\int^\infty_1 \tau^{-p_\theta-1}F(\tau)\,d\tau<\infty.
\end{equation}
\end{corollary}
{\bf Proof.}
Assume that problem~(P) possesses a solution with $\mu=\delta_y$ in $Q_T$
for some $T>0$.
Set
\begin{equation}
\label{eq:4.39}
f(\tau):=0\quad\mbox{for}\quad 0\le\tau\le R,
\quad
f(\tau):=F(\tau)-\tau^dR^{-d}F(R)\quad\mbox{for}\quad \tau>R.
\end{equation}
Then, by (D2)-(i) we see that $f$ is increasing and $F\ge f$ in $[0,\infty)$.
Applying Theorem~\ref{Theorem:3.1} 
with $z=y$, so that $m_\sigma(z)=\delta_y(B(y,\sigma ))\equiv 1$,  we find $\gamma\ge 1$ such that 
$$
\int_{\gamma^{-1}T^{-\frac{N}{\theta}}}^{\gamma^{-1}\sigma^{-N}}
s^{-p_\theta-1}f(s)\,ds\le \gamma^{p_\theta+1},
\qquad 0<\sigma<T^{\frac{1}{\theta}}.
$$
Letting $\sigma\to 0$, we have
$$
\int_{\gamma^{-1}T^{-\frac{N}{\theta}}}^\infty
s^{-p_\theta-1}f(s)\,ds\le \gamma^{p_\theta+1}.
$$
This together with \eqref{eq:4.39} implies \eqref{eq:4.38}.

Conversely, under condition~\eqref{eq:4.38},
we apply Theorem~\ref{Theorem:4.1} to obtain a local-in-time solution of problem~(P)
with $\mu=\delta_y$.
Thus Corollary~\ref{Corollary:4.4} follows.
\qed\vspace{5pt}

We mention that the integral condition \eqref{eq:4.38} also appears 
in \cite{LRSV}*{Theorem~5.1} as a necessary and sufficient conditions for existence with $L^1$ initial data.  
See also the informal argument preceding the proof of Theorem~4.1 of that work, where a Dirac delta function is considered as initial data. 
%
\section{Proof of the Main Theorem}
\label{section:5}
%
{\bf Proof of Theorem~\ref{Theorem:1.1}.}
Assertion~(i) is proved by Corollary~\ref{Corollary:3.1}~(ii), Remark~\ref{Remark:2.1}, and Corollary~\ref{Corollary:4.1}.

We now prove the non-existence parts of statements~(1) and~(2) in assertion~(ii). Suppose first that 
\eqref{eq:1.3} holds and there exists a local solution of problem (P).
Then 
$$
\sup_{z\in{{  \R}}^N}\mu(B(z,\sigma))
\ge \gamma_1\int_{B(0,\sigma)}|x|^{-N}|\log |x||^{-1}[\log|\log|x||]^{-\frac{N}{\theta}-1}\,dx
\ge C_1\gamma_1[\log|\log \sigma|]^{-\frac{N}{\theta}}$$
for small enough {$\sigma>0$}. 
For large enough $\gamma_1$ we then obtain a contradiction to Corollary~\ref{Corollary:3.1}. 
Hence no local solution can exist for such $\gamma_1$. Now suppose that \eqref{eq:1.5} holds. 
Then there exists $C_2>0$ such that 
$$
\sup_{z\in{{  \R}}^N}\mu(B(z,\sigma))
\ge \gamma_2\int_{B(0,\sigma)}|x|^{-N}|\log |x||^{-\frac{N(q+1)}{\theta}-1}\,dx
\ge C_2\gamma_2|\log \sigma|^{-\frac{N(q+1)}{\theta}}
$$
for small enough $\sigma>0$. 
Again we can obtain a contradiction to Corollary~\ref{Corollary:3.1} for large enough $\gamma_2$
and deduce that problem~{\rm (P)} possesses no local-in-time solution for such $\gamma_2$.  

Next we prove the existence parts of statements~(1) and (2) in assertion~(ii).  
Assume therefore that either \eqref{eq:1.4} with $\varepsilon_1\in(0,1)$ or \eqref{eq:1.6} with $\varepsilon_2\in(0,1)$ hold. 
Let $\alpha>0$ and set 
\begin{equation*}
h(\tau):=
\left\{
\begin{array}{ll}
\log(e+\log(e+\tau)) & \mbox{if \eqref{eq:1.4} holds},\vspace{5pt}\\ 
(\log(e+\tau))^{q+1} & \mbox{if \eqref{eq:1.6} holds},
\end{array}
\right.
\qquad
\psi_\alpha^\pm(\tau):=\tau h(\tau)^{\pm\alpha}.
\end{equation*}
If \eqref{eq:1.4} holds, then 
\begin{equation*}
\begin{split}
 \psi_\alpha^+(\mu(x)) & \le C\mu(x)\log(e+\log(e+\mu(x)))^\alpha\\
 & \le C\varepsilon_1|x|^{-N}|\log |x||^{-1}[\log|\log|x||]^{-\frac{N}{\theta}-1+\alpha}\chi_{B(0,R)}(x)+C,
 \quad x\in{  \R}^N.
\end{split}
\end{equation*}
This implies that 
\begin{equation*}
\begin{split}
\sup_{x\in{  \R}^N}\psi_\alpha^-\left[\,\dashint_{B(x,\sigma)}\psi_\alpha^+(\mu)\,dy\right]
 & \le \psi_\alpha^-(C\varepsilon_1\sigma^{-N}[\log|\log\sigma|]^{-\frac{N}{\theta}+\alpha})\\
 & \le C\varepsilon_1\sigma^{-N}[\log|\log\sigma|]^{-\frac{N}{\theta}}
\le C\varepsilon_1\sigma^{-N}h(\sigma^{-1})^{-\frac{N}{\theta}}
\end{split}
\end{equation*}
for small enough $\sigma>0$. 

Similarly, if \eqref{eq:1.6} holds, then 
\begin{equation*}
\begin{split}
 \psi_\alpha^+(\mu(x)) & \le C\mu(x)[\log(e+\mu(x))]^{\alpha(q+1)}\\
 & \le C\varepsilon_2|x|^{-N}|\log |x||^{-\frac{N(q+1)}{\theta}-1+\alpha(q+1)}\chi_{B(0,R)}(x)+C,
 \quad x\in{  \R}^N.
\end{split}
\end{equation*}
This implies that 
\begin{equation*}
\begin{split}
\sup_{x\in{  \R}^N}\psi_\alpha^-\left[\,\dashint_{B(x,\sigma)}\psi_\alpha^+(\mu)\,dy\right]
 & \le \psi_\alpha^-(C\varepsilon_2\sigma^{-N}|\log\sigma|^{-\frac{N(q+1)}{\theta}+\alpha(q+1)})\\
 & \le C\varepsilon_2\sigma^{-N}|\log\sigma|^{-\frac{N(q+1)}{\theta}}
\le C\varepsilon_2\sigma^{-N}h(\sigma^{-1})^{-\frac{N}{\theta}}
\end{split}
\end{equation*}
for small enough $\sigma>0$. 
Therefore, by Corollary~\ref{Corollary:4.2} we see that, 
if $\varepsilon_1>0$ (respectively $\varepsilon_2>0$) is small enough, 
then problem~{\rm (P)} possesses a local-in-time solution. 
Thus statements~(1) and~(2) in assertion~(ii) follow. 

Finally, we prove statement~(3) in assertion~(ii).
Assume that \eqref{eq:1.7} holds. 
Then 
$$
\sup_{z\in{{  \R}}^N}\mu(B(z,\sigma))
\ge \gamma_3\int_{B(0,\sigma)}|x|^{-\frac{\theta}{p-1}}|\log |x||^{-\frac{q}{p-1}}\,dx
\ge C_1\gamma_3\sigma^{N-\frac{\theta}{p-1}}|\log\sigma|^{-\frac{q}{p-1}}
$$
for small enough $\sigma>0$. 
This together with Corollary~\ref{Corollary:3.1} implies that 
problem~{\rm (P)} possesses no local-in-time solution for large enough $\gamma_3$. 
Conversely, suppose that \eqref{eq:1.8} holds. 
Since $p>p_\theta$, 
we find $\alpha>1$ such that $\alpha\theta/(p-1)<N$. 
Then we have 
$$
\sup_{x\in{{  \R}}^N}\left[\,\dashint_{B(x,\sigma)}
\mu(y)^\alpha\,dy\,\right]^{\frac{1}{\alpha}}
\le \left[C\varepsilon_3^\alpha\sigma^{-\frac{\alpha\theta}{p-1}}|\log\sigma|^{-\frac{\alpha q}{p-1}}+CK_3^\alpha\right]^{\frac{1}{\alpha}}
\le C\varepsilon_3\sigma^{-\frac{\theta}{p-1}}|\log\sigma|^{-\frac{q}{p-1}}
$$
for small enough $\sigma>0$. 
By Corollary~\ref{Corollary:4.3} we see that, 
if $\varepsilon_3>0$ is small enough, then problem~{\rm (P)} possesses a local-in-time solution. 
Thus statement~(3) in assertion~(ii) follows. 
The proof is complete.
\qed
\begin{remark}
\label{Remark:5.1}
The arguments in the proof of Theorem~{\rm\ref{Theorem:1.1}} are 
readily adapted to further log-refinements. For example, suppose that {\rm (F2)} is replaced by 
\begin{itemize}
  \item[{\rm(F2')}] 
  $F(\tau)\asymp \tau^p[\log\tau]^q[\log(\log\tau)]^r$ as $\tau\to\infty$ for some $p>1$ and $q$, $r\in{  \R}$. 
\end{itemize} 
Then we can show that problem~{\rm (P)} possesses a local-in-time solution if and only if 
$$
\sup_{z\in{  \R}^N}\mu(B(z,1))<\infty
$$ 
in the cases when {\rm (i)} $p<p_\theta$, {\rm (ii)} $p=p_\theta$ and $q<-1$, and {\rm (iii)} $p=p_\theta$, $q=-1$, and $r<-1$.
In the other cases, we divide condition~{\rm (F2')} into 4 cases: 
\begin{itemize}
  \item[{\rm (1)}]
  $p=p_\theta$ and $q=r=-1$.
  \item[{\rm (2)}]
  $p=p_\theta$, $q=-1$, and $r>-1$.
  \item[{\rm (3)}]
  $p=p_\theta$, $q>-1$, and $r\in{  \R}$.
  \item[{\rm (4)}]
  $p>p_\theta$ and $q$, $r\in{  \R}$,
\end{itemize}
and we can identify the optimal singularities of the initial data for solvability of problem~{\rm (P)}.  
Since inclusion of the proofs here would make the paper unduly long, 
 we leave the details to the interested reader.
\end{remark}
%

\noindent
{\bf Acknowledgment.}
YF was supported partially by JSPS KAKENHI Grant Number 19K14569. 
KH and KI  were supported in part by JSPS KAKENHI Grant Number JP19H05599. 
RL  was supported partially by a Daiwa Anglo-Japanese Foundation grant (Ref: 4646/13713).

\end{document}